\documentclass[10pt]{amsart}
\usepackage{amsmath}
\usepackage{amssymb}
\usepackage{enumerate}
\usepackage{amsbsy}
\usepackage{amsfonts}
\usepackage{color}

\newtheorem{thm}{Theorem}[section]
\newtheorem{lem}[thm]{Lemma}
\newtheorem{prop}[thm]{Propsition}

\newtheorem{defn}[thm]{Definition}
\newtheorem{rem}[thm]{Remark}

\numberwithin{equation}{section}

\begin{document}

	\title[Chern-Simons-Dirac System]{Well-Posedness in a critical space of Chern-Simons-Dirac System in the Lorenz gauge}

	\author[Y. Cho]{Yonggeun Cho}
	\address{Department of Mathematics, and Institute of Pure and Applied Mathematics, Jeonbuk National University, Jeonju 54896, Republic of Korea}
    \email{changocho@jbnu.ac.kr}

    \author[S. Hong]{Seokchang Hong}
    \address{Department of Mathematical Sciences, Seoul National University, Seoul 08826, Republic of Korea}
    \email{seokchangh11@snu.ac.kr}

		\begin{abstract}
		In this paper, we consider the Cauchy problem of local well-posedness of the Chern-Simons-Dirac system in the Lorenz gauge for $B^{\frac14}_{2,1}$ initial data. We improve the low regularity well-posedness, compared to Huh-Oh \cite{huhoh} and Okamoto \cite{oka}, by using the localization of space-time Fourier side and bilinear estimates given by Selberg \cite{selb}, whereas the authors of \cite{huhoh, oka} used global estimates of \cite{danfoselb}. Then we show the Dirac spinor flow of Chern-Simons-Dirac system is not $C^2$ at the origin in $H^s$ if $ s < \frac14$. From this point of view, the space $B_{2,1}^\frac14$ can be regarded as a critical space for the local well-posedness. We apply the argument for failure of smoothness to the Dirac equation decoupled from Chern-Simons-Dirac system and show the flow is not $C^3$ in $H^s, s < 0$.
	\end{abstract}

		\maketitle

\section{Introduction}
 We consider the Cauchy problem of Chern-Simons-Dirac (CSD) system with Lorenz gauge
 \begin{align}\label{csd}
 \left\{
 \begin{array}{l}
 \Box A_{\nu}  =  \partial^{\mu}(-2\epsilon_{\mu\nu\lambda}\psi^{\dagger}\alpha^{\lambda}\psi),\\
 (-i\alpha^{\mu}\partial_{\mu}+M\beta)\psi  =  A_{\mu}\alpha^{\mu}\psi, \\
 \psi(0) = \psi_0,\quad A_\mu(0) = a_\mu,\\
 \partial_0 A_0(0) = -\partial^ja_j, \quad \partial_0 A_j(0) = \partial_j a_0 - 2\epsilon_{0jk}\psi_0^\dagger\alpha^k\psi_0
 \end{array}\right.\end{align}
on the Minkowski space $\mathbb R^{1+2}$ equipped with the Minkowski metric of signature $(+, -, -)$. We adopt the Einstein summation convention, where Greek indices refer to $0, 1, 2$ and Latin indices refer to $1, 2$. Here $\Box = \partial^\mu \partial_\mu$ and $A_\mu: \mathbb R^{1+2} \to \mathbb R$ denotes the gauge potentials. The spinor field $\psi$ is represented by a column vector with two complex components. $\psi^\dagger$ is the complex conjugate transpose of $\psi$. $M \ge 0$ is the mass of the spinor field $\psi$, and $\epsilon_{\mu\nu\lambda}$ is the totally skew-symmetric tensor with $\epsilon_{012} = 1$. The Dirac matrices $\alpha$ and $\beta$ are given by
$$\alpha^{0} = I_{2\times2}, \quad \alpha^{1} = \sigma^{1}, \quad \alpha^{2} = \sigma^{2},\quad \beta = \sigma^{3},$$
where $\sigma^\mu$ are Pauli matrices such as
$$
\sigma^1 = \left(
          \begin{array}{cc}
            0 & 1 \\
            1 & 0 \\
          \end{array}
        \right), \quad \sigma^2 = \left(
          \begin{array}{cc}
            0 & -i \\
            i & 0 \\
          \end{array}
        \right),\quad \sigma^3 =
        \left(
          \begin{array}{cc}
            1 & 0 \\
            0 & -1 \\
          \end{array}
        \right).
$$

The (CSD) system \eqref{csd} is rewritten in the Lorenz gauge $\partial^\mu A_\mu = 0$ from the system of curvature form, which was introduced by Li-Bhaduri \cite{libha} and Cho, Kim, and Park \cite{ckp} to consider $(2+1)$-dimensional Chern-Simons electrodynamics (see also \cite{bm, dunne}), and it has the conservation of charge $Q(t) = \int_{\mathbb R^2} |\psi(t, x)|^2\,dx = Q(0)$
and $L^2$-scaling invariance in the case of $M = 0$.

In this paper we consider the local well-posedness (LWP) of \eqref{csd} in the Besov space $B_{2, 1}^\frac14$. The Besov space $B_{2, 1}^{\frac14}$ is defined by $\{f \in L^2 : \|f\|_{B_{2, 1}^\frac14} := \sum_{N :\, dyadic\; \ge 1}N^\frac14\|P_N f\|_{L_x^2} < \infty\}$, where $P_N$ is Littlewood-Paley projection on $\{|\xi| \sim N\}$. It is well-known that $H^s \subsetneq B_{2,1}^\frac14 \subsetneq H^\frac14$ for $s > \frac14$. The Sobolev space $H^s$ is defined by $\{f : \|f\|_{H^s} := (\sum N^{2s}\|P_Nf\|_{L_x^2}^2)^\frac12 < \infty\}$. %By Plancerel's theorem this Soblev norm is equivalent to $\|(1 + |\xi|)^s \widehat f \|_{L_{\xi}^2}$.
Our main result is stated as follows.
%%%%%%%%%%%%%%%%%%%%%%%%%%%%%%%%%%%%%%%%%%%%%%%%%%%%%%%%%%%%%%%%%%%%%%%%%%%%%%%%%%%%%%%%%%%%%%%%%%%%%%%%%%%%%%%%%%%%%%%%%%%%%%%%%%%%%%%%%%%%%%%%%%%%%%%%%%%%%%%%%%%%%%%%%%%%%%%%%%%%%%%%%%%%%%%%%%%%%%%%%%%%%%%%%%%%%%%%%%%%%%%%%%
%%%%%%%%%%%%%%%%%%%%%%%%%%%%%%%%%%%%%%%%%%%%%%%%%%%%%%%%%%%%%%%%%%%%%%%%%%%%%%%%%%%%%%%%%%%%%%%%%%%%%%%%%%%%%%%%%%%%%%%%%%%%%%%%%%%%%%%%%%%%%%%%%%%%%%%%%%%%%%%%%%%%%%%%%%%%%%%%%%%%%%%%%%%%%%%%%%%%%%%%%%%%%%%%%%%%%%%%%%%%%%%%%%
\begin{thm}\label{lwp}
Suppose that $a_0, a_1, a_2, \psi_0 \in B_{2, 1}^\frac14$. Then there exists $T = T(\|a_\mu\|_{B_{2, 1}^\frac14}, \|\psi_0\|_{B_{2, 1}^\frac14}, M) > 0$ such that there exist unique solution $A_\mu, \psi \in C((-T,T); B_{2,1}^\frac14)$ of \eqref{csd}, which depends continuously on the initial data.
\end{thm}
%%%%%%%%%%%%%%%%%%%%%%%%%%%%%%%%%%%%%%%%%%%%%%%%%%%%%%%%%%%%%%%%%%%%%%%%%%%%%%%%%%%%%%%%%%%%%%%%%%%%%%%%%%%%%%%%%%%%%%%%%%%%%%%%%%%%%%%%%%%%%%%%%%%%%%%%%%%%%%%%%%%%%%%%%%%%%%%%%%%%%%%%%%%%%%%%%%%%%%%%%%%%%%%%%%%%%%%%%%%%%%%%%%
%%%%%%%%%%%%%%%%%%%%%%%%%%%%%%%%%%%%%%%%%%%%%%%%%%%%%%%%%%%%%%%%%%%%%%%%%%%%%%%%%%%%%%%%%%%%%%%%%%%%%%%%%%%%%%%%%%%%%%%%%%%%%%%%%%%%%%%%%%%%%%%%%%%%%%%%%%%%%%%%%%%%%%%%%%%%%%%%%%%%%%%%%%%%%%%%%%%%%%%%%%%%%%%%%%%%%%%%%%%%%%%%%%

The Cauchy problem for \eqref{csd} was studied for the local well-posedness (LWP) in $H^s$ for $s > \frac14$ by Huh-Oh \cite{huhoh} and Okamoto \cite{oka}, and recently in Fourier-Lebesgue space $\widehat{H}^{s, r}$ for $1 < r \le 2$ and $s > \frac3{2r} - \frac12$ by Pecher \cite{pech}. Since $H^s \subsetneq B_{2,1}^\frac14$ for $s > \frac14$, Theorem \ref{lwp} improves the previous regularity and uniqueness results. For more related results we refer the readers to the references in \cite{huhoh, oka, pech}.

For the proof of Theorem \ref{lwp} we develop a duality argument in Besov type $X^{s, b}$ space (see \eqref{norm-rec} and Sections 6, 7 below). While doing this we must encounter the serious trilinear terms. In order to handle trilinear terms we reveal the null structures introduced in \cite{seltes, huhoh} and exploit the angular Whitney decomposition and $2D$ bilinear estimates of Selberg \cite{selb}. In the previous works \cite{huhoh, oka} the authors used the global bilinear estimates of \cite{danfoselb} or partial local estimates of \cite{danselb} and square-sum $(\ell^2)$ on the frequencies, which results in the well-posedness in $H^s, s > \frac14$. To improve the regularity and uniqueness we use the modulation and frequency localization in space-time Fourier side fully. To control the low-low-high modulation case we find an exclusion condition among frequencies (see Remark \ref{llh-mod}) and apply the angular Whitney decomposition. For the summation on the modulation we need to use $\ell^1$ summation instead of $\ell^2$, which gives us the well-posedness in $B_{2, 1}^\frac14$ instead of $H^\frac14$. Such Besov type well-posedness can be also found in Dirac-Klein-Gordon system \cite{wang}.

We now make a remark on the smoothness of the flow map $(a_\nu,\psi_0) \mapsto (A_\nu(t), \psi(t))$. The nonlinear term of \eqref{csd} is essentially quadratic. Hence one may expect the flow will be smooth in local time if the problem is well-posed. However, such smoothness can be shown to fail when the initial data are rougher than $H^\frac14$. For the simplicity sake we only consider the massless case $M = 0$. The massive case $M > 0$ can also be treated similarly with the phase $e^{\pm it\sqrt{M^2 + |\xi|^2}}$.
\begin{thm}\label{fail}
Let $M = 0$, $s < \frac14$ and $T > 0$. Then the flow map of \eqref{csd} $(a_\nu, \psi_0) \mapsto (A_\nu, \psi)$ from $H^s(\mathbb R^2)$ to $C([-T, T]; H^s(\mathbb R^2))$ cannot be $C^2$ at the origin.
\end{thm}
In view of Theorem \ref{fail} the Sobolev index $\frac14$ is critical. We can say that $B_{2, 1}^\frac14$ is a critical space and our regularity result is very sharp. It would be very interesting to find out a new method to resolve the regularity problem in $H^\frac14$. In \cite{oka} the author considered the same problem and showed that the flow of vector potential $A$ is not $C^2$. In this paper, we show the smoothness failure of  the flow for Dirac spinor $\psi$. We can take advantage of the degree of freedom for initial data with respect to not only $\psi$ and but also $A$. We proceed the proof by following the well-known argument in \cite{mst} and \cite{hele}.

On the other hand, one may consider the single equation of Dirac spinor field $\psi$ decoupled from \eqref{csd} with $M = 0$ and $a_\mu = 0$:
\begin{align}\label{single}
-i\alpha^{\nu}\partial_{\nu}\psi  =  \Box^{-1} \left[\partial^{\mu}(-2\epsilon_{\mu\nu\lambda}\psi^{\dagger}\alpha^{\lambda}\psi)\right]\alpha^{\nu}\psi,\quad \psi(0) = \psi_0 \in H^s.
\end{align}
Indeed, the author of \cite{oka} considered this type equation for the local well-posedness (see also \cite{danselb}). We show that the flow is not $C^3$ as above if the data are rougher than the scaling critical space $L^2$.
\begin{thm}\label{illp}
Let $s < 0$ and $T > 0$. Then the flow map of \eqref{single} $\psi_0 \mapsto \psi$ from $H^s(\mathbb R^2)$ to $C([-T, T]; H^s(\mathbb R^2))$ cannot be $C^3$ at the origin.
\end{thm}
The main difficulty is coming from the upper bound of phase multiplier
$$
\mathbf m_{1\cdots 5} = \int_0^t \!\!\!\int_0^{t'} e^{-\pm_1i(t-t')|\xi| -\pm_2i(t'-t'')|\eta|+ \pm_3it''|\zeta-\eta| - \pm_4it''|\zeta| - \pm_5it'|\xi-\eta|}\,dt''dt'.
$$
The essential part is the case $\pm_1 = \pm_5 , \pm_3 = \pm_4$, with $\pm_1 \neq \pm_3$, which results in $|\mathbf m_{1=5,3=4}| \lesssim |\eta|$. This makes it hard to  control the lower bound for second derivative of flow and consequently yields the failure when $s < 0$.

%On the other hand, if we consider the stationary field $A$, that is, $\partial_t A = 0$, then

Our paper is organized as follows. In Section 2, we give some preliminaries on Dirac projection operator, decomposition of d'Alemvertian, and Besov type $X^{s,b}$ space. In section 3, we introduce bilinear estimates based on the null structure. In Section 4, we give a sketch of proof of Theorem \ref{lwp} by constructing Picard's iterates. In Sections 5,6,7, we give proofs of crucial parts for LWP. In Section 8, we prove the failure of smoothness \eqref{single} for the system when $s < \frac14$ and for the decoupled equation when $s < 0$. The last section is devoted to proving the energy estimate, Lemma \ref{energy-est}.\\

%%%%%%%%%%%%%%%%%%%%%%%%%%%%%%%%%%%%%%%%%%%%%%%%%%%%%%%%%%%%%%%%%%%%%%%%%%%%%%%%%%%%%%%%%%%%%%%%%%%%%%%%%%%%%%%%%%%%%%%%%%%%%%%%%%%%%%%%%%%%%%%%%%%%%%%%%%%%%%%%%%%%%%%%%%%%%%%%%%%%%%%%%%%%%%%%%%%%%%%%%%%%%%%%%%%%%%%%%%%%%%%%%%%%%%%%%%%%%%%%%%%%%%%%%%%%%%%%%%%%%%%%%%%%%%%%%%%%%%%%%%%%%%%%%%%%%%%%%%%%%%%%%%%%%%%%%%%%%%%%%%%%%%%%%%%%%%%%%%%%%%%%%%%%%%%%%%%%%%%%%%%%%%%%%%%%%%%%%%%%%%%%%%%%%%%%%%%%%%%%%%%%%%%%%%%%%%%%%%%%%%%%%%%%%%%%%%%%%%%%%%%%%%%%%%%%%%%%%%%%%%%%%%%%%%%%%%%%%%%%%%%%%%%%%%%%%%%%%%%%%%%%%%%%%%%%%%%%%%%%%%%%%%%%%%%%%%%%%%%%%%%%%%%%%%%%%%%%%%%%%%%%%%%%%%%%%%%%%%%%%%%%%%%%%%%%%%%%%%%%%%%%%%%%%%%%%%%%%%%%%%%%%%%%%%%%
\noindent\textbf{Notations.}\\
%%%%%%%%%%%%%%%%%%%%%%%%%%%%%%%%%%%%%%%%%%%%%%%%%%%%%%%%%%%%%%%%%%%%%%%%%%%%%%%%%%%%%%%%%%%%%%%%%%%%%%%%%%%%%%%%%%%%%%%%
\noindent$\bullet$ Since we often use $L_{t, x}^{2}$ norm, we abbreviate $\|F\|_{L^{2}_{t,x}}$ by $\|F\|$.\\
%%%%%%%%%%%%%%%%%%%%%%%%%%%%%%%%%%%%%%%%%%%%%%%%%%%%%%%%%%%%%%%%%%%%%%%%%%%%%%%%%%%%%%%%%%%%%%%%%%%%%%%%%%%%%%%%%%%%%%%%
\noindent$\bullet$ The spatial Fourier transform and space-time Fourier transform on $\mathbb{R}^{2}$ and $\mathbb R^{1+2}$ are defined by
$$\widehat{f}(\xi)=\int_{\mathbb{R}^{2}}e^{-ix\cdot\xi}f(x)dx,\quad \widetilde{u}(X)=\int_{\mathbb R^{1+2}}e^{-i(t\tau+x\cdot\xi)}u(t,x)dtdx,$$
where $\tau\in\mathbb{R}$, $\xi\in\mathbb{R}^{2}$, and $X = (\tau,\xi)\in\mathbb R^{1+2}$. Also we denote $\mathcal{F}(u)=\widetilde{u}$.\\
%%%%%%%%%%%%%%%%%%%%%%%%%%%%%%%%%%%%%%%%%%%%%%%%%%%%%%%%%%%%%%%%%%%%%%%%%%%%%%%%%%%%%%%%%%%%%%%%%%%%%%%%%%%%%%%%%%%%%%%%%
\noindent$\bullet$ We denote $D := |\nabla|$ whose symbol is $|\xi|$.\\
%%%%%%%%%%%%%%%%%%%%%%%%%%%%%%%%%%%%%%%%%%%%%%%%%%%%%%%%%%%%%%%%%%%%%%%%%%%%%%%%%%%%%%%%%%%%%%%%%%%%%%%%%%%%%%%%%%%%%%%%%
\noindent$\bullet$ For any $E \subset \mathbb R^{1+2}$ the projection operator $P_E$ is defined by $\widetilde{P_E u}(\tau, \xi) = \chi_E \widetilde u(\tau, \xi)$.\\
%%%%%%%%%%%%%%%%%%%%%%%%%%%%%%%%%%%%%%%%%%%%%%%%%%%%%%%%%%%%%%%%%%%%%%%%%%%%%%%%%%%%%%%%%%%%%%%%%%%%%%%%%%%%%%%%%%%%%%%%%
\noindent$\bullet$ As usual different positive constants depending only on $M$ are denoted by the same letter $C$, if not specified. $A \lesssim B$ and $A \gtrsim B$ means that $A \le CB$ and
$A \ge C^{-1}B$, respectively for some $C>0$. $A \sim B$ means that $A \lesssim B$ and $A \gtrsim B$.\\

\section{Preliminaries}

\subsection{Dirac projection operator}
We first define Dirac projection operator by
$$\Pi_{\pm}(\xi)=\frac{1}{2}\left(I_{2\times 2} \pm \frac{\xi_{j}\alpha^{j}}{|\xi|}\right).$$
Then we have useful identity
$$\alpha^{i}\Pi_{\pm}=\Pi_{\mp}\alpha^{i}+\frac{\xi^{i}}{|\xi|} I_{2 \times 2}.$$
Also we define Riesz transform $R_{\pm}^{\mu}$, $\mu=0,1,2$ by
$$R^{0}_{\pm}=-1, \quad R_{\pm}^{j}=\mp\frac{\partial_{j}}{iD},\quad R_{\pm, j} := - R_{\pm}^j.$$
Using the above identities, we get
\begin{align}\label{id-proj}
\alpha^{\mu}\Pi_{\pm}=\Pi_{\mp}\alpha^{\mu}\Pi_{\pm}-R^{\mu}_{\pm}\Pi_{\pm}.
\end{align}

Let $\psi_\pm := \Pi_{\pm}\psi$ and $\psi_{\pm}^{\rm hom} = e^{\mp itD}\psi_{0}$. Then Duhamel's principle gives us the following integral equation:
\begin{align}\label{psi}
\psi_{\pm}(t,x)  =  \psi_{\pm}^{\rm hom} - i\int_{0}^{t}e^{\mp i(t-t')D}M\beta\psi_{\mp}(t',x)dt'+i\int_{0}^{t}e^{\mp i(t-t')D}\Pi_{\pm}(A_{\mu}\alpha^{\mu}\psi)(t',x)dt'.
\end{align}
Since we have from \eqref{id-proj} that
\begin{eqnarray*}
\Pi_{\pm}(A_{\mu}\alpha^{\mu}\Pi_{\pm}\psi_{\pm}) & = & \Pi_{\pm}(A_{\mu}(\Pi_{\mp}\alpha^{\mu}\Pi_{\pm}-R_{\pm}^{\mu}\Pi_{\pm})\psi_{\pm}) \\
& = & \Pi_{\pm}A_{\mu}\Pi_{\mp}\alpha^{\mu}\psi_{\pm}-\Pi_{\pm}A_{\mu}R_{\pm}^{\mu}\psi_{\pm},
\end{eqnarray*}
we write
\begin{align}\begin{aligned}\label{dirac-eqn}
\psi_{\pm}(t,x) = \psi_{\pm}^{\rm hom}  &-  i\int_{0}^{t}e^{\mp i(t-t')D}M\beta\psi_{\mp}(t',x)dt' \\
&+  i\int_{0}^{t}e^{\mp i(t-t')D}\Pi_{\pm}(A_{\mu}\Pi_{\mp}\alpha^{\mu}\Pi_{\pm}\psi_{\pm})(t',x)dt' \\
&-  i\int_{0}^{t}e^{\mp i(t-t')D}\Pi_{\pm}A_{\mu}R_{\pm}^{\mu}\psi_{\pm}(t',x)dt'.
\end{aligned}\end{align}

\subsection{Decomposition of d'Alembertian}

\begin{lem}[Lemma 2.1 of \cite{huhoh}]\label{decom-box}
Suppose that $\Box u = \partial^{\mu}F_{\mu}$ for Schwartz functions $u, F_\mu$. Then $u = u_+ + u_-$, where
$$
u_{\pm}(t,x) = u_\pm^{\rm hom}(t, x) - \frac{1}{2}\int_{0}^{t}e^{\mp i(t-t')D}R_{\pm}^{\mu}F_{\mu}(t',x)dt',
$$
where
$$
u_\pm^{\rm hom}(t, x) = \frac12\left( e^{\mp itD}u(0) \pm \frac1{iD}(F_0(0) - \partial_0u(0))\right).
$$
\end{lem}
From wave part of \eqref{csd} and Lemma \ref{decom-box} we may write
\begin{eqnarray}
A_{\nu,\pm}(t,x)  =  A_{\nu,\pm}^{\rm hom}(t,x) - \frac{1}{2}\int_{0}^{t}e^{\mp i(t-t')D}R_{\pm}^{\mu}(-2\epsilon_{\mu\nu\lambda}\psi^{\dagger}\alpha^{\lambda}\psi)(t',x)dt',\label{anu}
\end{eqnarray}
where
$$
A_{\nu,\pm}^{\rm hom}(t,x) = \frac12 e^{\mp itD}\left(A_\nu(0) \pm \frac1{iD}(-\epsilon_{0\nu\lambda}\psi(0)^\dagger \alpha^\lambda \psi(0) - \partial_0 A_\nu(0))\right).
$$
%where we assume that
%$$\psi_{0}\in B^{1/4}_{2,1}(\mathbf{R}^{2},\mathbf{C}^{2}), \quad  a_{\mu}(x)=A_{\mu}(0,x)\in B^{1/4}_{2,1}(\mathbf{R}^{2},\mathbf{R}^{2})$$
Now we exploit the nonlinear terms:
\begin{eqnarray*}
\epsilon_{\mu\nu\lambda}\psi_{1,\pm_{1}}^{\dagger}\alpha^{\lambda}\psi_{2,\pm_{2}} & = & \epsilon_{\mu\nu\lambda}\psi_{1,\pm_{1}}^{\dagger}\alpha^{\lambda}\Pi_{\pm_{2}}\psi_{2,\pm_{2}} \\
& = & \epsilon_{\mu\nu\lambda}\psi_{1,\pm_{1}}^{\dagger}(\Pi_{\mp_{2}}\alpha^{\lambda}\Pi_{\pm_{2}}-R_{\pm_{2}}^{\lambda}\Pi_{\pm_{2}})\psi_{2,\pm_{2}} \\
& = &  \epsilon_{\mu\nu\lambda}\psi_{1,\pm_{1}}^{\dagger}\Pi_{\mp_{2}}\alpha^{\lambda}\Pi_{\pm_{2}}\psi_{2,\pm_{2}}-\epsilon_{\mu\nu\lambda}\psi_{1,\pm_{1}}^{\dagger}R_{\pm_{2}}^{\lambda}\psi_{2,\pm_{2}}.
\end{eqnarray*}
Thus we write
\begin{align}\begin{aligned}\label{wave-eqn}
A_{\nu,\pm}(t,x) = A_{\nu,\pm}^{\rm hom}(t,x) &+ \sum_{\pm_1,\pm_2}\int_{0}^{t}e^{\mp i(t-t')D}R_{\pm}^{\mu}\epsilon_{\mu\nu\lambda}\psi_{1,\pm_{1}}^{\dagger}\Pi_{\mp_{2}}\alpha^{\lambda}\Pi_{\pm_{2}}\psi_{2,\pm_{2}}(t',x)dt' \\
&- \sum_{\pm_1,\pm_2}\int_{0}^{t}e^{\mp i(t-t')D}R_{\pm}^{\mu}\epsilon_{\mu\nu\lambda}\psi_{1,\pm_{1}}^{\dagger}R_{\pm_{2}}^{\lambda}\psi_{2,\pm_{2}}(t',x)dt'.
\end{aligned}\end{align}

\subsection{Function spaces}

Now we define function spaces.
\begin{defn}
Let $s, b \in \mathbb R$. Let $N,L\ge 1$ be dyadic number. We define Besov type $X^{s, b}$ space by
$$\|u\|_{\mathcal{B}^{s,b;1}_{\pm}}=\sum_{N,L\ge1}N^{s}L^{b}\|P_{K_{N,L}^\pm} u\|$$ and
$$\|u\|_{\mathcal{B}^{s,b;\infty}_{\pm}}=\sup_{N,L\ge1}N^{s}L^{b}\|P_{K_{N,L}^\pm}u\|, $$
where $$K_{N,L}^{\pm}= \{ (\tau,\xi)\in\mathbb R^{1+2} : |\xi|\sim N,\;\; |\tau + \pm |\xi| |\sim L \}.$$
\end{defn}
Then the $\mathcal{B}^{s,b;1}_{\pm}$ norm can be recovered as follows:
\begin{equation}
\|u\|_{\mathcal{B}^{s,b;1}_{\pm}} = \sup_{\|v\|_{\mathcal{B}^{-s,-b;\infty}_{\pm}}=1}\bigg|\int u\bar{v}dtdx\bigg|.\label{norm-rec}
\end{equation}

Now we consider the time-slab
$$
S_{T}=(-T,T)\times\mathbb{R}^{2}.
$$
Define
$$
\|u\|_{\mathcal{B}^{s, b; 1}_{\pm}(S_{T})}=\inf_{v = u \text{ on }S_{T}}\|v\|_{\mathcal{B}^{s, b; 1}_{\pm}}.
$$
This becomes a semi-norm on $\mathcal{B}^{s, b; 1}_{\pm}$, but is a norm if we identify elements which agree on $S_{T}$ and the resulting space is denoted $\mathcal{B}^{s, b; 1}_{\pm}(S_{T})$. In other words, $\mathcal{B}^{s, b; 1}_{\pm}(S_{T})$ is the quotient space $\mathcal{B}^{s,b;p}_{\pm}/X$, where $X= \{ v\in\mathcal{B}^{s,b;p}_{\pm } : v = 0 \text{ on } S_{T} \}$. Since $X$ is a closed subspace in $\mathcal{B}^{s, b; 1}_{\pm}$, we conclude that the quotient $\mathcal{B}^{s, b; 1}_{\pm}(S_{T})$ is a Banach space.

To control the nonlinear terms of \eqref{dirac-eqn} and \eqref{wave-eqn} in $\mathcal{B}^{s, b; 1}_{\pm}(S_{T})$, we introduce the energy estimate lemma.
\begin{lem}\label{energy-est}
Let us consider the integral equation:
$$
v(t) = e^{\mp itD}f + \int_{0}^{t}e^{\mp i(t-t')D}F(t')dt'.
$$
with sufficiently smooth $f$ and $F$.
If $T \le 1$, then for any $s \in \mathbb R$ we have
$$
\|v\|_{\mathcal{B}^{s,\frac{1}{2};1}_{\pm}(S_{T})} \lesssim  \|f\|_{B^{s}_{2,1}} + \|F\|_{\mathcal{B}^{s,-\frac{1}{2};1}_{\pm}(S_{T})}.
$$
\end{lem}
 The proof is exactly same as the one of \cite{danselb}. But for the convenience of readers, we append its proof in the last section.

%%%%%%%%%%%%%%%%%%%%%%%%%%%%%%%%%%%%%%%%%%%%%%%%%%%%%%%%%%%%%%%%%%%%%%%%%%%%%%%%%%%%%%%%%%%%%%%%%%%%%%%%%%%%%%%%%%%%%%%%%%%%%%%%%%%%%%%%%%%%%%%%%%%%%%%%%%%%%%%%%

%%%%%%%%%%%%%%%%%%%%%%%%%%%%%%%%%%%%%%%%%%%%%%%%%%%%%%%%%%%%%%%%%%%%%%%%%%%%%%%%%%%%%%%%%%%%%%%%%%%%%%%%%%%%%%%%%%%%%%%%%%%%%%%%%%%%%%%%%%%%%%%%%%%%%%%%%%%%%%%%%

\section{Bilinear estimates and null structure}

\subsection{Bilinear estimates}
For dyadic $N,L\ge 1$, let us invoke that
$$
K_{N,L}^{\pm} = \{ (\tau,\xi)\in\mathbb R^{1+2} : |\xi|\sim N,\;\; |\tau \pm |\xi| |\sim L \}.
$$

To handle the nonlinear terms in \eqref{dirac-eqn} and \eqref{wave-eqn}, we utilize the $2$-dimensional bilinear estimates of wave type shown by Selberg.
\begin{thm}[Theorem 2.1 of \cite{selb}]\label{bilinear-est}
For all $u_{1},u_{2}\in L^{2}_{t,x}(\mathbb{R}^{1+2})$ such that $\widetilde{u_{j}}$ is supported in $K_{N_{j},L_{j}}^{\pm_{j}}$, the estimate
$$
\|P_{K_{N_{0},L_{0}}^{\pm_{0}}}(u_{1}\overline{u_{2}})\| \le C\|u_{1}\|\|u_{2}\|
$$
holds with
\begin{eqnarray}
C & \sim & (N_{\min}^{012}L_{\min}^{12})^{1/2}(N_{\min}^{12}L_{\max}^{12})^{1/4},\label{bi-selb-1} \\
C & \sim & (N_{\min}^{012}L_{\min}^{0j})^{1/2}(N_{\min}^{0j}L_{\max}^{0j})^{1/4}, \quad j = 1, 2,\label{bi-selb-2} \\
C & \sim & ((N_{\min}^{012})^{2}L_{\min}^{012})^{1/2}\label{bi-selb-3}
\end{eqnarray}
regardless of the choices of signs $\pm_{j}$.
\end{thm}
In this paper, we only use \eqref{bi-selb-1} and \eqref{bi-selb-2} to show Theorem \ref{lwp}. In dealing with these estimates, we assume $L^{12}_{\max}\ll N^{012}_{\min}$. Otherwise, the estimate \eqref{bi-selb-3} is better than \eqref{bi-selb-1} and \eqref{bi-selb-2}.

In the proof of Theorem \ref{lwp}, we exclude the low-low-high modulation with high-high-low frequency and high-low-high frequency.(See Remark \ref{llh-mod}.) Then we are left to deal with the low-low-high modulation with all input and output frequencies compatible. For this purpose, we first apply the angular Whitney decomposition of \cite{selb} as follows:
For $\gamma, r > 0$ and $\omega\in\mathbb S^1$, where $\mathbb S^1 \subset \mathbb R^2$ is the unit circle, we define
\begin{align*}
\Gamma_\gamma(\omega)&= \{ \xi\in\mathbb R^2 : \angle(\xi,\omega)\le\gamma \},\\
T_r(\omega) &= \{ \xi\in\mathbb R^2 : |P_{\omega^\perp}\xi|\lesssim r \},
\end{align*}
where $P_{\omega^\perp}$ is the projection onto the orthogonal complement $\omega^\perp$ of $\omega$ in $\mathbb R^2$. Also we let $\Omega(\gamma)$ denote a maximal $\gamma$-separated subset of the unit circle. Then %for $\xi_1,\xi_2\in\mathbb R^2\setminus\{0\}$ with $\angle(\omega_1,\omega_2)>0$, we have
%\begin{align}\label{ang-whi1}
%1\sim \sum_{\substack{ 0<\gamma<1 \\ \gamma:{\rm dyadic} }}\sum_{\substack{ \omega_1,\omega_2\in\Omega(\gamma) \\ 3\gamma\le\angle(\omega_1,\omega_2)\le12\gamma }}\chi_{\Gamma_\gamma(\omega_1)}(\xi_1)\chi_{\Gamma_\gamma(\omega_2)}(\xi_2).
%\end{align}
%Note that the lower bound $\angle(\omega_1,\omega_2)\ge3\gamma$ says that the sectors $\Gamma_\gamma(\omega_1)$ and $\Gamma_\gamma(\omega_2)$ are well-separated. If this separation is not required, we then have the following modification:
for $0<\gamma<1$ and $k\in\mathbb N$, we have
\begin{align}\label{ang-whi2}
\chi_{\angle(\xi_1,\xi_2)\le k\gamma} \lesssim \sum_{\substack{ \omega_1,\omega_2\in\Omega(\gamma) \\ \angle(\omega_1,\omega_2)\le(k+2)\gamma }}\chi_{\Gamma_\gamma(\omega_1)}(\xi_1)\chi_{\Gamma_\gamma(\omega_2)}(\xi_2),
\end{align}
for all $\xi_1,\xi_2\in\mathbb R^2\setminus\{0\}$ with $\angle(\omega_1,\omega_2) > 0$.

Second, we apply the following null form estimate.
\begin{thm}[Theorem 2.3 of \cite{selb}]\label{null-form-selb}
Let $r>0$ and $\omega\in\mathbb S^1$. Then for all $u_1,u_2\in L^2_{t,x}(\mathbb R^{1+2})$ such that $\widetilde{u_j}$ is supported in $K_{N_j,L_j}^{\pm_j}$, we have
$$
\| B_{\theta_{12}}(P_{T_r(\omega)}u_1,u_2)\| \lesssim (rL_1L_2)^{1/2}\|u_1\|\|u_2\|.
$$
\end{thm}
Here, the bilinear form $ B_{\theta_{12}}(u_1,u_2)$ is defined on the Fourier side by inserting the angle $\theta_{12}=\angle(\pm_1\xi_1,\pm_2\xi_2)$ in the convolution of $u_1$ and $u_2$; that is,
$$
\mathcal F B_{\theta_{12}}(u_1,u_2)(X_0) = \int_{X_0=X_1+X_2}\angle(\pm_1\xi_1,\pm_2\xi_2)\widetilde{u_1}(X_1)\widetilde{u_2}(X_2)dX_1dX_2.
$$

\subsection{Bilinear interaction}
The space-time Fourier transform of the product $\psi_{2}^{\dagger}\psi_{1}$ of two spinor fields $\psi_{1}$ and $\psi_{2}$ is written as
$$
\widetilde{\psi_{2}^{\dagger}\psi_{1}}(X_{0})=\int_{X_{0} = X_{1} - X_{2}}\widetilde{\psi_{2}}^{\dagger}(X_{2})\widetilde{\psi_{1}}(X_{1})\, dX_{1}dX_{2},
$$
where $\psi^{\dagger}$ is the transpose of complex conjugate of $\psi$. Here the relation between $X_{1}$ and $X_{2}$ in the convolution integral of spinor fields is given by $X_{0} = X_{1} - X_{2}$ so called bilinear interaction. This is also the case for the product of two complex scalar fields. %Motivated by this relation in convolution integral, we call a triple $(X_{0},X_{1},X_{2})$ of vectors $X_{j}=(\tau_{j},\xi_{j})\in\mathbb{R}^{1+2}$ to be a bilinear interaction if $X_{0} = X_{1}-X_{2}$ holds. Given signs $(\pm_{0},\pm_{1},\pm_{2})$ we also define the hyperbolic weights $h_{j}=\tau_{j}\pm_{j}|\xi_{j}|$. If all $h_{j}$ vanish, we say that the interaction is null. In this case, all $X_{j}$ lie on the null cone $\{ (\tau,\xi)\in\mathbb{R}^{1+2} : \tau = \pm|\xi| \}$ and $\angle(\pm_{1}\xi_{1},\pm_{2}\xi_{2})=0$.
The following lemma is on the bilinear interaction.
\begin{lem}[Lemma 2.2 of \cite{selb}]\label{bi-int}
Given a bilinear interaction $(X_{0},X_{1},X_{2})$ with $\xi_{j} \neq 0$, and signs $(\pm_{0},\pm_{1},\pm_{2})$, let $h_{j}=\tau_{j}\pm_{j}|\xi_{j}|$ and $\theta_{12} = |\angle(\pm_{1}\xi_{1},\pm_{2}\xi_{2})|$. Then we have
\begin{equation}
\max(|h_{0}|,|h_{1}|,|h_{2}|) \gtrsim \min(|\xi_{1}|,|\xi_{2}|)\theta_{12}^{2}.
\end{equation}
\end{lem}
Note that using interpolation with (3.4) and the trivial inequality $\theta_{12}\le 1$, one can obtain
\begin{equation}
\theta_{12}\lesssim \bigg(\frac{\max(|h_{0}|,|h_{1}|,|h_{2}|)}{\min(|\xi_{1}|,|\xi_{2}|)}\bigg)^{p}, \quad 0\le p\le\frac{1}{2}.\label{bi-int-int}
\end{equation}

\subsection{Null structure}
Let us consider bilinear forms defined by
\begin{eqnarray*}
Q^{\mu\nu}_{\pm_{1},\pm_{2}}(\phi_{1,\pm_{1}},\phi_{2,\pm_{2}}) & = & R^{\mu}_{\pm_{1}}\phi_{1,\pm_{1}}R^{\nu}_{\pm_{2}}\phi_{2,\pm_{2}}-R^{\nu}_{\pm_{1}}\phi_{1,\pm_{1}}R^{\mu}_{\pm_{2}}\phi_{2,\pm_{2}}, \\
Q^{0}_{\pm_{1},\pm_{2}}(\phi_{1,\pm_{1}},\phi_{2,\pm_{2}}) & = & R^{\mu}_{\pm_{1}}\phi_{1,\pm_{1}} R_{\mu,\pm_{2}}\phi_{2,\pm_{2}}.
\end{eqnarray*}
 \begin{lem}[Lemma 2.6 of \cite{huhoh}]
 Let $\phi_{1},\phi_{2}$ be complex-valued Schwartz functions, and $\psi_{1},\psi_{2}$ Schwartz spinor fields. Then we have
 \begin{eqnarray}
 |\mathcal{F}\lbrack Q^{\mu\nu}_{\pm_{1},\pm_{2}}(\phi_{1},\phi_{2})\rbrack(X_{0})| & \lesssim & \int_{X_{0} = X_{1} + X_{2}}\theta_{12}\,|\widetilde{\phi_{1}}(X_{1})||\widetilde{\phi_{2}}(X_{2})|\,dX_{1}dX_{2},\label{bil-1} \\
|\mathcal{F}\lbrack Q^{0}_{\pm_{1},\pm_{2}}(\phi_{1},\phi_{2})\rbrack(X_{0})| & \lesssim & \int_{X_{0} = X_{1} + X_{2}} \theta_{12}^{2} |\widetilde{\phi_{1}}(X_{1})||\widetilde{\phi_{2}}(X_{2})|\,dX_{1}dX_{2}, \label{bil-2} \\
|\mathcal{F}\lbrack(\Pi_{\pm_{1}}\psi_{1})^{\dagger}(\Pi_{\mp_{2}}\alpha^{\mu}\Pi_{\pm_{2}}\psi_{2})(X_{0})| & \lesssim & \int_{X_{0} = X_{1} + X_{2}}\theta_{12}\,|\widetilde{\psi_{1}}(X_{1})||\widetilde{\psi_{2}}(X_{2})|\,dX_{1}dX_{2}.\label{bil-3}
 \end{eqnarray}
 \end{lem}

\

%%%%%%%%%%%%%%%%%%%%%%%%%%%%%%%%%%%%%%%%%%%%%%%%%%%%%%%%%%%%%%%%%%%%%%%%%%%%%%%%%%%%%%%%%%%%%%%%%%%%%%%%%%%%%%%%%%%%%%%%%%%%%%%%%%%%%%%%%%%%%%%%%%%%%%%%%%%%%%%%%%%%%%%%%%%%%%%%%%%%%%%%%%%%%%%%%%%%%%%%%%%%%%%%%%%%%%%%%%%%%%%%%%%%%%%%%%%%%%%%%%%%%%%%%%%%%%%%%%%%%%%%%%%%%%%%%%%%%%%%%%%%%%%%%%%%%%%%%%%%%%%%%%%%%%%%%%%%%%%%%%%%%%%%%%%%%%%%%%%%%%%%%%%%%%%%%%%%%%%%%%%%%%%%%%%%%%%%%%%%%%%%%%%%%%%%%%%%%%%%%%%%%%%%%%%%%%%%%%%%%%%%%%%%%%%%%%%%%%%%%%%%%%%%%%%%%%%%%%%%%%%%%%%%%%%%%%%%%%%%%%%%%%%%%%%%%%%%%%%%%%%%%%%%%%%%%%%%%%%%%%%%%%%%%%%%%%%%%%%%%%%%%%%%%%%%%%%%%%%%%%%%%%%%%%%%%%%%%%%%%%%%%%%%%%%%%%%%%%%%%%%%%%%%%%%%%%%%%%%%%%%%%%%%%%%%%%%%%%%%%%%%%%%%%%%%%%%%%%%%%%%%%%%%%%%%%%%%%%%%%%%%%%%%%%%%%%%%%%%%%%%%%%%%%%%%%%%%%%%%%%%%%%%%%%%%%%%%%%%%%%%%%%%%%%%%%%%%%%%%%%%%%%%%
\section{Sketch of proof of Theorem \ref{lwp}}
We first construct a Picard's iterates for the equations \eqref{dirac-eqn} and \eqref{wave-eqn} as follows.
Set $A_{\nu,\pm}^{(0)} = A_{\nu,\pm}^{\rm hom}, \psi_{\pm}^{(0)} = \psi_{\pm}^{\rm hom}$, and for $n\ge1$ define $A_{\nu,\pm}^{(n)},\psi_{\pm}^{(n)}$ as
\begin{eqnarray*}
A^{(n)}_{\nu,\pm}(t,x) & = & A^{(0)}_{\nu,\pm}(t,x) + \int_{0}^{t}e^{\mp(t-t')D}R_{\pm}^{\mu}\mathcal{N}_{\mu\nu}(\psi^{(n-1)},\psi^{(n-1)})\,dt', \\
\psi^{(n)}_{\pm}(t,x) & = & \psi^{(0)}_{\pm}(t,x) - i\int_{0}^{t}e^{\mp(t-t')D}\left(M\beta\psi^{(n-1)}_\mp-\Pi_{\pm}\mathcal{M}(\psi^{(n-1)},A^{(n-1)})\right)\,dt',
\end{eqnarray*}
 where
 \begin{align*}
 &\mathcal{N}_{\mu\nu}(\psi_{1,\pm_1},\psi_{2,\pm_2}) = \epsilon_{\mu\nu\lambda}\psi_{1,\pm_{1}}^{\dagger}\Pi_{\mp_{2}}\alpha^{\lambda}\Pi_{\pm_{2}}\psi_{2,\pm_{2}} - \epsilon_{\mu\nu\lambda}\psi_{1,\pm_{1}}^{\dagger}R_{\pm_{2}}^{\lambda}\psi_{2,\pm_{2}},\\
 & \mathcal{M}(\psi,A) = -A_{\mu}\Pi_{\mp}\alpha^{\mu}\Pi_{\pm}\psi_{\pm} + A_{\mu}R_{\pm}^{\mu}\psi_{\pm}.
 \end{align*}

Then we show that Picard's iteration converges. For this purpose we need to prove that $(A_{\nu,+}^{(n)},A_{\nu,-}^{(n)},\psi^{(n)}_{+},\psi^{(n)}_{-})$ is a Cauchy sequence in the space $\mathcal{B}^{\frac{1}{4},\frac{1}{2};1}_{+}\times\mathcal{B}^{\frac{1}{4},\frac{1}{2};1}_{-}\times\mathcal{B}^{\frac{1}{4},\frac{1}{2};1}_{+}\times\mathcal{B}^{\frac{1}{4},\frac{1}{2};1}_{-}$.
In fact, we have only to show that the following estimates:
\begin{align}
&\|A^{\rm hom}_{\nu,\pm}\|_{\mathcal{B}^{\frac{1}{4},\frac{1}{2};1}_{\pm}(S_{T})} \lesssim \sum_{\mu=0,1,2}\|a_{\mu}\|_{B^{\frac14}_{2,1}}, \quad \|\psi^{\rm hom}_{\pm}\|_{\mathcal{B}^{\frac{1}{4},\frac{1}{2};1}_{\pm}(S_{T})} \lesssim \|\psi_{0}\|_{B^{\frac14}_{2,1}},\label{hom}\\
&\|M\beta\psi_{\mp}\|_{\mathcal{B}^{\frac{1}{4},-\frac{1}{2};1}_{\mp}(S_T)}  \lesssim  M\|\psi_{\mp}\|_{\mathcal{B}^{\frac{1}{4},\frac{1}{2};1}_{\mp}},\label{lin}\\
&\|R^{\mu}_{\pm}\mathcal{N}_{\mu\nu}(\psi_{1,\pm_{1}},\psi_{2,\pm_{2}})\|_{\mathcal{B}^{\frac{1}{4},-\frac{1}{2};1}_{\pm}(S_T)}  \lesssim  \|\psi_{1,\pm_{1}}\|_{\mathcal{B}^{\frac{1}{4},\frac{1}{2};1}_{\pm_{1}}}\|\psi_{2,\pm_{2}}\|_{\mathcal{B}^{\frac{1}{4},\frac{1}{2};1}_{\pm_{2}}}, \label{psipsi}\\
&\|\Pi_{\pm}\mathcal{M}(\psi_{1,\pm_{1}},A_{\mu,\pm_{2}})\|_{\mathcal{B}^{\frac{1}{4},-\frac{1}{2};1}_{\pm}(S_T)}  \lesssim
\|\psi_{1,\pm_{1}}\|_{\mathcal{B}^{\frac{1}{4},\frac{1}{2};1}_{\pm_{1}}}\|A_{\mu,\pm_{2}}\|_{\mathcal{B}^{\frac{1}{4},\frac{1}{2};1}_{\pm_{2}}}. \label{psia}
\end{align}
Now the proof of local well-posedness in $B^{\frac14}_{2,1}$ is quite standard. Thus we omit it. In the next three sections we will focus on the estimates above.

%%%%%%%%%%%%%%%%%%%%%%%%%%%%%%%%%%%%%%%%%%%%%%%%%%%%%%%%%%%%%%%%%%%%%%%%%%%%%%%%%%%%%%%%%%%%%%%%%%%%%%%%%%%%%%%%%%%%%%%%%%%%%%%%%%%%%%%%%%%%%%%%%%%%%%%%%%%%%%%%%%%%%%%%%%%%%%%%%%%%%%%%%%%%%%%%%%%%%%%%%%%%%%%%%%%%%%%%%%%%%%
%%%%%%%%%%%%%%%%%%%%%%%%%%%%%%%%%%%%%%%%%%%%%%%%%%%%%%%%%%%%%%%%%%%%%%%%%%%%%%%%%%%%%%%%%%%%%%%%%%%%%%%%%%%%%%%%%%%%%%%%%%%%%%%%%%%%%%%%%%%%%%%%%%%%%%%%%%%%%%%%%%%%%%%%%%%%%%%%%%%%%%%%%%%%%%%%%%%%%%%%%%%%%%%%%%%%%%%%%%%%%%
\section{Proof of \eqref{hom} and \eqref{lin}}

The initial condition of the (CSD) system \eqref{csd} says that
$$
\partial_{0}A(0,x) = -\partial^{j}a_{j}(x), \quad \partial_{t}A_{j}(0,x) = \partial_{j}a_{0}(x) - 2\epsilon_{0jk}\psi_{0}^{\dagger}\alpha^{k}\psi_{0}.
$$
Using these, $A^{\rm hom}_{\nu,\pm}$ can be rewritten as
\begin{eqnarray*}
A^{\rm hom}_{0,\pm}(t,x) & = & \frac{1}{2}e^{\mp itD}\bigg(a_{0}\pm\frac{\partial^{j}}{iD}a_{j}\bigg), \\
A^{\rm hom}_{j,\pm}(t,x) & = & \frac{1}{2}e^{\mp itD}\bigg(a_{j}\mp \frac{\partial_{j}}{iD}a_{0} \bigg).
\end{eqnarray*}
Now by Lemma \ref{energy-est}, we have
$$\|A^{\rm hom}_{0,\pm}\|_{\mathcal{B}^{\frac{1}{4},\frac{1}{2};1}_{\pm}(S_{T})} + \|A^{\rm hom}_{j,\pm}\|_{\mathcal{B}^{\frac{1}{4},\frac{1}{2};1}_{\pm}(S_{T})} \lesssim  \bigg\|a_{0}\pm\frac{\partial^{j}}{iD}a_{j}\bigg\|_{B^{\frac14}_{2,1}} + \bigg\|a_{j}\pm\frac{\partial^{j}}{iD}a_{0}\bigg\|_{B^{\frac14}_{2,1}}
 \lesssim  \sum_{\mu=0,1,2}\|a_{\mu}\|_{B^{\frac14}_{2,1}},$$
where we used the fact that the Riesz transform $\partial_{j}/iD$ is bounded on $L^{2}_{x}$. $\psi^{\rm hom}_{\pm}$ can be treated similarly, and the estimate \eqref{lin} for linear term follows immediately from the definition of $\mathcal{B}^{s, b ; 1}_{\pm}$ and their embedding.

%%%%%%%%%%%%%%%%%%%%%%%%%%%%%%%%%%%%%%%%%%%%%%%%%%%%%%%%%%%%%%%%%%%%%%%%%%%%%%%%%%%%%%%%%%%%%%%%%%%%%%%%%%%%%%%%%%%%%%%%%%%%%%%%%%%%%%%%%%%%%%%%%%%%%%%%%%%%%%%%%%%%%%%%%%%%%%%%%%%%%%%%%%%%%%%%%%%%%%%%%%%%%%%%%%%%%%%%%%%%%%
%%%%%%%%%%%%%%%%%%%%%%%%%%%%%%%%%%%%%%%%%%%%%%%%%%%%%%%%%%%%%%%%%%%%%%%%%%%%%%%%%%%%%%%%%%%%%%%%%%%%%%%%%%%%%%%%%%%%%%%%%%%%%%%%%%%%%%%%%%%%%%%%%%%%%%%%%%%%%%%%%%%%%%%%%%%%%%%%%%%%%%%%%%%%%%%%%%%%%%%%%%%%%%%%%%%%%%%%%%%%%%

\section{Estimates of $\mathcal N_{\mu\nu}$: Proof of \eqref{psipsi}}

In this section, we prove \eqref{psipsi}. To this end we need to show that
\begin{eqnarray}
\|R_{\pm}^{\mu}\epsilon_{\mu\nu\lambda}\psi_{1,\pm_{1}}^{\dagger}\Pi_{\mp_{2}}\alpha^{\lambda}\Pi_{\pm_{2}}\psi_{2,\pm_{2}}\|_{\mathcal{B}^{\frac{1}{4},-\frac{1}{2};1}_{\pm}}  & \lesssim & \prod_{j=1,2}\|\psi_{j,\pm_{j}}\|_{\mathcal{B}^{\frac{1}{4},\frac{1}{2};1}_{\pm_{j}}},\label{nl-anu-1}  \\
\|R_{\pm}^{\mu}\epsilon_{\mu\nu\lambda}\psi_{1,\pm_{1}}R^{\lambda}_{\pm_{2}}\psi_{2,\pm_{2}}\|_{\mathcal{B}^{\frac{1}{4},-\frac{1}{2};1}_{\pm}}  & \lesssim & \prod_{j=1,2}\|\psi_{j,\pm_{j}}\|_{\mathcal{B}^{\frac{1}{4},\frac{1}{2};1}_{\pm_{j}}}.\label{nl-anu-2}
\end{eqnarray}

\subsection{Proof of \eqref{nl-anu-1}}
The duality formula \eqref{norm-rec} yields that
$$
\|R_{\pm}^{\mu}\epsilon_{\mu\nu\lambda}\psi_{1,\pm_{1}}^{\dagger}\Pi_{\mp_{2}}\alpha^{\lambda}\Pi_{\pm_{2}}\psi_{2,\pm_{2}}\|_{\mathcal{B}^{\frac{1}{4},-\frac{1}{2};1}_{\pm}}
 = \sup_{\|\phi\|_{\mathcal{B}^{-\frac{1}{4},\frac{1}{2};\infty}_{\pm}} = 1}\left|\int_{\mathbb R^{1+2}}  R_{\pm}^{\mu}\epsilon_{\mu\nu\lambda}\psi_{1,\pm_{1}}^{\dagger}\Pi_{\mp_{2}}\alpha^{\lambda}\Pi_{\pm_{2}}\psi_{2,\pm_{2}}\overline{\phi}\,dtdx\right|.
$$
By self-adjointness of Riesz transform we have
\begin{align*}
\int  R_{\pm}^{\mu}\epsilon_{\mu\nu\lambda}\psi_{1,\pm_{1}}^{\dagger}\Pi_{\mp_{2}}\alpha^{\lambda}\Pi_{\pm_{2}}\psi_{2,\pm_{2}}\overline{\phi}\,dtdx  &= \int  \epsilon_{\mu\nu\lambda}\psi_{1,\pm_{1}}^{\dagger}\Pi_{\mp_{2}}\alpha^{\lambda}\Pi_{\pm_{2}}\psi_{2,\pm_{2}}\overline{R^{\mu}_{\pm}\phi}\,dtdx \\
&= \int \overline{\mathcal F(R^{\mu}_{\pm}\phi)}(-X_0)\,\epsilon_{\mu\nu\lambda}\mathcal F(\psi_{1,\pm_{1}}^{\dagger}\Pi_{\mp_{2}}\alpha^{\lambda}\Pi_{\pm_{2}}\psi_{2,\pm_{2}}) (X_0) \,dX_{0} =: \mathbf J^1.
\end{align*}
A dyadic decomposition of space-time Fourier side gives us $|\mathbf J^1| \le \sum_{\mathbf{N,L}} |\mathbf J^1_{\mathbf{N,L}}|$, where
$$
|\mathbf J^1_{\mathbf{N,L}}| \le  |\epsilon_{\mu\nu\lambda}|\left|\int \overline{\mathcal F(R^{\mu}_{\pm} P_{K_{N_{0},L_{0}}^{\pm}}\phi)}(-X_0)\,\mathcal F((P_{K_{N_{1},L_{1}}^{\pm_{1}}}\psi_{1,\pm_{1}})^{\dagger}\Pi_{\mp_{2}}\alpha^{\lambda}\Pi_{\pm_{2}}P_{K_{N_{2},L_{2}}^{\pm_{2}}}\psi_{2,\pm_{2}})(X_0) \,dX_{0}\right|,
$$
and
$$
\mathbf{N}=(N_{0},N_{1},N_{2}), \quad \mathbf{L}=(L_{0},L_{1},L_{2}).
$$
We can also write $|\mathbf J^1_{\mathbf{N,L}}|$ as
$$
|\epsilon_{\mu\nu\lambda}|\left|\int \mathcal F(P_{K_{N_{1},L_{1}}^{\pm_{1}}}\psi_{1,\pm_{1}})^{\dagger}(-X_0)
\,\mathcal F( \overline{R^{\mu}_{\pm} P_{K_{N_{0},L_{0}}^{\pm}}\phi}\,\Pi_{\mp_{2}}\alpha^{\lambda}\Pi_{\pm_{2}}P_{K_{N_{2},L_{2}}^{\pm_{2}}}\psi_{2,\pm_{2}})(X_0) \,dX_{0}\right|,
$$

\begin{rem}\label{llh-mod}
Here, we make a remark on a low-low-high modulation and low-high-high frequency case $L_{{\rm max}}^{12} \ll L_0 \ll N_1 \ll N_0 \sim N_2$. In view of the first representation of $\mathbf J^1_{\mathbf{N,L}}$, the support condition excludes the case $\pm = \pm_2$ since $\phi$ and $\psi_{2, \pm_2}$ have up and down cones or down and up cones, respectively. On the other hand, the support condition of the second representation excludes the case $\pm \neq \pm_2$ since $\phi$ and $\psi_{2, \pm_2}$ have now up and up, or down and down cones, respectively. Therefore we can conclude that this type llh-lhh case does not appear in the summation of $\mathbf J^1_{\mathbf{N,L}}$. By the same way, we can exclude the case llh-hhl $L_{{\rm max}}^{12} \ll L_0 \ll N_0 \ll N_1 \sim N_2$.
\end{rem}

%the choices of signs $\pm_j$ and $\pm_j'$, $j=0,1,2$. By definition of $\mathcal B^{s,b,;1}_\pm$, we must have $\pm=\pm_0'$ in $\mathbf J^1$. On the other hand, we see that
%\begin{align*}
%\|\psi_+\|_{\mathcal B^{s,b,;1}_-} &= \sum_{N,L\ge1}N^sL^b\|P_{K_{N,L}^-}\Pi_+\psi\| \\
%&= \sum_{N,L\ge1}N^sL^b\|\chi_{K_{N,L}^-}(\tau,\xi)\Pi_+(\xi)\widetilde{\psi}(\tau,\xi)\| \\
%&= \sum_{N,L\ge1}N^sL^b\|\chi_{|\xi|\sim N}(\xi)\chi_{|\tau-|\xi||\sim L}(\tau,\xi)(I_{2\times2}+\frac{\xi_j\alpha^j}{|\xi|})\widetilde{\psi}(\tau,\xi)\|
%\end{align*}

We can write the integrand of $\mathbf J^1$ with the combination of positive and negative parts of every component of $\mathcal F P_{K_{N_{0},L_{0}}^{\pm}}\phi$, $\mathcal F P_{K_{N_{1},L_{1}}^{\pm_{1}}}\psi_{1,\pm_{1}}$, and $\mathcal F P_{K_{N_{2},L_{2}}^{\pm_{2}}}\psi_{2,\pm_{2}}$. Thus it is harmless to assume that $\mathcal F P_{K_{N_{0},L_{0}}^{\pm}}\phi$, $\mathcal F P_{K_{N_{1},L_{1}}^{\pm_{1}}}\psi_{1,\pm_{1}}$, and $\mathcal F P_{K_{N_{2},L_{2}}^{\pm_{2}}}\psi_{2,\pm_{2}}$ are nonnegative real-valued functions.

Since $R^{0}_{\pm}=-1$ and $|\widetilde{R^{j}_{\pm}\phi}|\le |\widetilde{\phi}|$,  using the bilinear estimate \eqref{bil-3}, we get
$$
|\mathbf J^1| \lesssim \sum_{\mathbf{N,L}}\int \mathcal F P_{K_{N_{0},L_{0}}^{\pm}}\phi \int_{X_{0}=X_{1}-X_{2}}\theta_{12}\;\chi_{K_{N_{0},L_{0}}^{\pm}}(X_1-X_2) \mathcal F P_{K_{N_{1},L_{1}}^{\pm_{1}}}\psi_{1,\pm_{1}}(X_{1})\,\mathcal F P_{K_{N_{2},L_{2}}^{\pm_{2}}}\psi_{2,\pm_{2}}(X_{2})\, dX_{1}dX_{2}dX_{0},
$$
By Cauchy-Schwarz inequality and Lemma \ref{bi-int},
$$
|\mathbf J^1| \lesssim  \sum_{\mathbf{N,L}}\|P_{K_{N_{0},L_{0}}^{\pm}}\phi\|\left(\frac{L_{\max}^{012}}{N_{\min}^{12}}\right)^{1/2}
\|P_{K_{N_{0},L_{0}}^{\pm}}(P_{K_{N_{1},L_{1}}^{\pm_{1}}}\psi_{1,\pm_{1}}P_{K_{N_{2},L_{2}}^{\pm_{2}}}\psi_{2,\pm_{2}})\|.
$$

\subsubsection{Case 1: $L_{0}\le L_{1}\le L_{2}$}\hfill

If $N_{0}\ll N_{1}\sim N_{2}$, then we use \eqref{bi-selb-2} with  $j = 1$  of Theorem \ref{bilinear-est} and get
\begin{eqnarray*}
|\mathbf J^1| & \lesssim & \sum_{\mathbf{N,L}}\|P_{K_{N_{0},L_{0}}^{\pm}}\phi\|\bigg(\frac{L_{2}}{N_{1}}\bigg)^{1/2}(N_{0}L_{0})^{1/2}(N_{0}L_{1})^{1/4}\|P_{K_{N_{1},L_{1}}^{\pm_{1}}}\psi_{1,\pm_{1}}\|\|P_{K_{N_{2},L_{2}}^{\pm_{2}}}\psi_{2,\pm_{2}}\| \\
& = & \sum_{\mathbf{N,L}}\|P_{K_{N_{0},L_{0}}^{\pm}}\phi\|\,N_{0}^{1/4}L_{0}^{-1/2}\bigg(\frac{L_{2}}{N_{1}}\bigg)^{1/2}(N_{0}L_{0})^{1/2}(N_{0}L_{1})^{1/4} \|P_{K_{N_{1},L_{1}}^{\pm_{1}}}\psi_{1,\pm_{1}}\|\|P_{K_{N_{2},L_{2}}^{\pm_{2}}}\psi_{2,\pm_{2}}\| \\
& \lesssim & \|\phi\|_{\mathcal{B}^{-\frac{1}{4},\frac{1}{2};\infty}_{\pm}}\sum_{\mathbf{N,L}}\bigg(\frac{N_{0}}{N_{1}}\bigg)^{1/2}N_{1}^{1/4}N_{2}^{1/4}L_{1}^{1/4}L_{2}^{1/2} \|P_{K_{N_{1},L_{1}}^{\pm_{1}}}\psi_{1,\pm_{1}}\|\|P_{K_{N_{2},L_{2}}^{\pm_{2}}}\psi_{2,\pm_{2}}\|  \\
& \lesssim & \|\psi_{1,\pm_{1}}\|_{\mathcal{B}^{\frac{1}{4},\frac{1}{2};1}_{\pm_{1}}}\|\psi_{2,\pm_{2}}\|_{\mathcal{B}^{\frac{1}{4},\frac{1}{2};1}_{\pm_{2}}}\|\phi\|_{\mathcal{B}^{-\frac{1}{4},\frac{1}{2};\infty}_{\pm}}.
\end{eqnarray*}
Here, we used $$\sum_{L_{0};L_{0}\le L_{1}}\sum_{N_{0};N_{0}\ll N_{1}}\bigg(\frac{N_{0}}{N_{1}}\bigg)^{1/2}\lesssim L_{1}^{\epsilon}.$$

If $N_{1}\lesssim N_{0}\sim N_{2}$, then by \eqref{bi-selb-2} we also get
\begin{align*}
|\mathbf J^1| &\lesssim \|\phi\|_{\mathcal{B}^{-\frac{1}{4},\frac{1}{2};\infty}_{\pm}}\sum_{\mathbf{N,L}}N_{0}^{1/4}L_{0}^{-1/2}\bigg(\frac{L_{2}}{N_{1}}\bigg)^{1/2}(N_{1}L_{0})^{1/2}(N_{1}L_{1})^{1/4} \|P_{K_{N_{1},L_{1}}^{\pm_{1}}}\psi_{1,\pm_{1}}\|\|P_{K_{N_{2},L_{2}}^{\pm_{2}}}\psi_{2,\pm_{2}}\| \\
&= \|\phi\|_{\mathcal{B}^{-\frac{1}{4},\frac{1}{2};\infty}_{\pm}}\sum_{\mathbf{N,L}}N_{0}^{1/4}L_{2}^{1/2}N_{1}^{1/4}L_{1}^{1/4} \|P_{K_{N_{1},L_{1}}^{\pm_{1}}}\psi_{1,\pm_{1}}\|\|P_{K_{N_{2},L_{2}}^{\pm_{2}}}\psi_{2,\pm_{2}}\| \\
& \lesssim  \|\psi_{1,\pm_{1}}\|_{\mathcal{B}^{\frac{1}{4},\frac{1}{2};1}_{\pm_{1}}}\|\psi_{2,\pm_{2}}\|_{\mathcal{B}^{\frac{1}{4},\frac{1}{2};1}_{\pm_{2}}}\|\phi\|_{\mathcal{B}^{-\frac{1}{4},\frac{1}{2};\infty}_{\pm}} \;\;\left(\because \sum_{L_{0};L_{0}\le L_{1}}\sum_{N_{0};N_{0}\sim N_{2}}N_{0}^\frac14 \lesssim  N_{2}^{1/4}L_{1}^{\epsilon}\right)	.
\end{align*}

\subsubsection{Case 2: $L_{1}\le L_{0}\le L_{2}$}\hfill

For $N_{0}\ll N_{1}\sim N_{2}$ it follows from  \eqref{bi-selb-2} that
\begin{align*}
|\mathbf J^1| &\lesssim \|\phi\|_{\mathcal{B}^{-\frac{1}{4},\frac{1}{2};\infty}_{\pm}}\sum_{\mathbf{N,L}}N_{0}^{1/4}L_{0}^{-1/2}\bigg(\frac{L_{2}}{N_{1}}\bigg)^{1/2}(N_{0}L_{1})^{1/2}(N_{0}L_{0})^{1/4} \|P_{K_{N_{1},L_{1}}^{\pm_{1}}}\psi_{1,\pm_{1}}\|\|P_{K_{N_{2},L_{2}}^{\pm_{2}}}\psi_{2,\pm_{2}}\| \\
&\lesssim  \|\phi\|_{\mathcal{B}^{-\frac{1}{4},\frac{1}{2};\infty}_{\pm}}\sum_{\mathbf{N,L}}\bigg(\frac{N_{0}}{N_{1}}\bigg)^{1/2}L_{0}^{-1/4}N_{1}^{1/4}N_{2}^{1/4}L_{1}^{1/2}L_{2}^{1/2} \|P_{K_{N_{1},L_{1}}^{\pm_{1}}}\psi_{1,\pm_{1}}\|\|P_{K_{N_{2},L_{2}}^{\pm_{2}}}\psi_{2,\pm_{2}}\| \\
&\lesssim  \|\psi_{1,\pm_{1}}\|_{\mathcal{B}^{\frac{1}{4},\frac{1}{2};1}_{\pm_{1}}}\|\psi_{2,\pm_{2}}\|_{\mathcal{B}^{\frac{1}{4},\frac{1}{2};1}_{\pm_{2}}}\|\phi\|_{\mathcal{B}^{-\frac{1}{4},\frac{1}{2};\infty}_{\pm}}  \;\;\left(\because \sum_{L_{0}\ge1}\sum_{N_{0};N_{0}\ll N_{1}}\bigg(\frac{N_{0}}{N_{1}}\bigg)^{1/2}L_{0}^{-1/4}\lesssim1\right)	.
\end{align*}
If $N_{1}\lesssim N_{0}\sim N_{2}$, then similarly
\begin{align*}
|\mathbf J^1| &\lesssim  \|\phi\|_{\mathcal{B}^{-\frac{1}{4},\frac{1}{2};\infty}_{\pm}}\sum_{\mathbf{N,L}}N_{0}^{1/4}L_{0}^{-1/2}\bigg(\frac{L_{2}}{N_{1}}\bigg)^{1/2}(N_{1}L_{1})^{1/2}(N_{1}L_{0})^{1/4}\|P_{K_{N_{1},L_{1}}^{\pm_{1}}}\psi_{1,\pm_{1}}\|\|P_{K_{N_{2},L_{2}}^{\pm_{2}}}\psi_{2,\pm_{2}}\| \\
&= \|\phi\|_{\mathcal{B}^{-\frac{1}{4},\frac{1}{2};\infty}_{\pm}} \sum_{\mathbf{N,L}}N_{0}^{1/4}L_{0}^{-1/4}L_{2}^{1/2}L_{1}^{1/2}N_{1}^{1/4} \|P_{K_{N_{1},L_{1}}^{\pm_{1}}}\psi_{1,\pm_{1}}\|\|P_{K_{N_{2},L_{2}}^{\pm_{2}}}\psi_{2,\pm_{2}}\|\\
& \lesssim  \|\psi_{1,\pm_{1}}\|_{\mathcal{B}^{\frac{1}{4},\frac{1}{2};1}_{\pm_{1}}}\|\psi_{2,\pm_{2}}\|_{\mathcal{B}^{\frac{1}{4},\frac{1}{2};1}_{\pm_{2}}}\|\phi\|_{\mathcal{B}^{-\frac{1}{4},\frac{1}{2};\infty}_{\pm}}	 \;\;\left(\because \sum_{L_{0}\ge1}\sum_{N_{0};N_{0}\sim N_{2}}N_{0}^{1/4}L_{0}^{-1/4}\lesssim N_{2}^{1/4}\right).
\end{align*}

\subsubsection{Case 3: $L_{1}\le L_{2}\le L_{0}$}\hfill

 %Now we consider the support condition carefully. First, if $L_2\ll L_0$, the most serious case is that $L_2\ll L_0\ll N_{\min}^{012}$. But by suitable choices of signs $\pm_j'$ we see that the integral vanishes in view of support condition. Indeed, if $N_0\ll N_1\sim N_2$, then by choosing $\pm_0'\neq (\pm_1'=\pm_2')$, the supports are disjoint. Also if $N_1\ll N_0\sim N_2$, then the sign $\pm_0'\neq \pm_2'$ says that the supports are disjoint. For $N_0\sim N_1\sim N_2$, with $\pm_0'\neq (\pm_1'=\pm_2')$ we see that the integral vanishes.

First, we consider the low-low-high modulation (i.e., $L_{\max}^{12}\ll L_0 \ll N_{{\rm min}}^{012}$). In view of Remark \ref{llh-mod} we have only to consider the case $N_0 \sim N_1 \sim N_2$. Let us invoke the angular Whitney decompositions \eqref{ang-whi2} to get
$$
|\mathbf J^1| \lesssim \sum_{\mathbf{N,L}}\sum_{\omega_1,\omega_2}\int \theta_{12}\mathcal F [(P_{K_{N_1,L_1}^{\pm_1}}\psi_{1,\pm_1}^{2\theta_{12},\omega_1})^\dagger P_{K_{N_2,L_2}^{\pm_2}}\psi_{2,\pm_2}^{2\theta_{12},\omega_2}](X_0)\overline{\mathcal F [R^\mu_{\pm}P_{K_{N_0,L_0}^\pm}\phi]}(-X_0)dX_0,
$$
where $\psi_{j,\pm_j}^{\theta,\omega_j}=P_{\pm_j\xi_j\in\Gamma_{\theta}(\omega_j)}\psi_{j,\pm_j}$. We note that the spatial Fourier support of $P_{K_{N_1,L_1}^{\pm_1}}\psi_{1,\pm_1}^{2\theta_{12},\omega_1}$ is contained in a strip of radius compatible to $N_{\max}^{12}\theta_{12}$ around $\mathbb R\omega_1$.

Then we use Theorem \ref{null-form-selb} with $r\sim N_{\max}^{12}\theta_{12} $ followed by Cauchy-Schwarz inequality to obtain
\begin{align*}
|\mathbf J^1| &\lesssim \sum_{\mathbf{N,L}}\sum_{\omega_1,\omega_2}\|P_{K_{N_0,L_0}^\pm}B_{\theta_{12}}(P_{K_{N_1,L_1}^{\pm_1}}\psi_{1,\pm_1}^{2\theta_{12},\omega_1},P_{K_{N_2,L_2}^{\pm_2}}\psi_{2,\pm_2}^{2\theta_{12},\omega_2})\|\|P_{K_{N_0,L_0}^\pm}\phi\| \\
& \lesssim \|\phi\|_{\mathcal{B}^{-\frac{1}{4},\frac{1}{2};\infty}_{\pm}}\sum_{\mathbf{N,L}}\sum_{\omega_1,\omega_2}N_0^{1/4}L_0^{-1/2}\left(N_{\max}^{12}\left(\frac{L_0}{N_{\min}^{12}}\right)^{1/2}L_1L_2\right)^{1/2}\|P_{K_{N_1,L_1}^{\pm_1}}\psi_{1,\pm_1}^{2\theta_{12},\omega_1}\|\|P_{K_{N_2,L_2}^{\pm_2}}\psi_{2,\pm_2}^{2\theta_{12},\omega_2}\| \\
& \lesssim \|\phi\|_{\mathcal{B}^{-\frac{1}{4},\frac{1}{2};\infty}_{\pm}}\sum_{\mathbf{N,L}}N_0^{1/2}L_0^{-1/4}(L_1L_2)^{1/2}\|P_{K_{N_1,L_1}^{\pm_1}}\psi_{1,\pm_1}\|\|P_{K_{N_2,L_2}^{\pm_2}}\psi_{2,\pm_2}\|,
\end{align*}
where the summation by $\omega_1,\omega_2\in\Omega(\theta_{12})$ is taken over $\angle(\omega_1,\omega_2)\le4\theta_{12}$.

Second, we consider the case $L_0 \gtrsim N_{\min}^{012}$. If $N_0\ll N_1\sim N_2$, then
\begin{align*}
|\mathbf J^1| &\lesssim  \|\phi\|_{\mathcal{B}^{-\frac{1}{4},\frac{1}{2};\infty}_{\pm}}\sum_{\mathbf{N,L}}N_{0}^{1/4}L_{0}^{-1/2}\left(\frac{L_{0}}{N_{1}}\right)^{1/2}(N_{0}L_{1})^{1/2}(N_{1}L_{2})^{1/4}\|P_{K_{N_{1},L_{1}}^{\pm_{1}}}\psi_{1,\pm_{1}}\|\|P_{K_{N_{2},L_{2}}^{\pm_{2}}}\psi_{2,\pm_{2}}\| 	\\
& \lesssim \|\phi\|_{\mathcal{B}^{-\frac{1}{4},\frac{1}{2};\infty}_{\pm}}\sum_{\mathbf{N,L}}N_{0}^{1/4}N_{0}^{-1/2}\left(\frac{L_{0}}{N_{1}}\right)^{1/2}(N_{0}L_{1})^{1/2}(N_{1}L_{2})^{1/4}\|P_{K_{N_{1},L_{1}}^{\pm_{1}}}\psi_{1,\pm_{1}}\|\|P_{K_{N_{2},L_{2}}^{\pm_{2}}}\psi_{2,\pm_{2}}\| \\
&= \|\phi\|_{\mathcal{B}^{-\frac{1}{4},\frac{1}{2};\infty}_{\pm}}\sum_{\mathbf{N,L}}N_{0}^{1/4}\left(\frac{L_{0}}{N_{1}}\right)^{1/2}L_{1}^{1/2}(N_{1}L_{2})^{1/4}\|P_{K_{N_{1},L_{1}}^{\pm_{1}}}\psi_{1,\pm_{1}}\|\|P_{K_{N_{2},L_{2}}^{\pm_{2}}}\psi_{2,\pm_{2}}\| \\
& \lesssim \|\psi_{1,\pm_{1}}\|_{\mathcal{B}^{\frac{1}{4},\frac{1}{2};1}_{\pm_{1}}}\|\psi_{2,\pm_{2}}\|_{\mathcal{B}^{\frac{1}{4},\frac{1}{2};1}_{\pm_{2}}}\|\phi\|_{\mathcal{B}^{-\frac{1}{4},\frac{1}{2};\infty}_{\pm}}.
\end{align*}

Also, for $N_1\lesssim N_0\sim N_2$, we see that
\begin{align*}
|\mathbf J^1| &\lesssim  \|\phi\|_{\mathcal{B}^{-\frac{1}{4},\frac{1}{2};\infty}_{\pm}}\sum_{\mathbf{N,L}}N_{0}^{1/4}L_{0}^{-1/2}\left(\frac{L_{0}}{N_{1}}\right)^{1/2}(N_{1}L_{1})^{1/2}(N_{1}L_{2})^{1/4}\|P_{K_{N_{1},L_{1}}^{\pm_{1}}}\psi_{1,\pm_{1}}\|\|P_{K_{N_{2},L_{2}}^{\pm_{2}}}\psi_{2,\pm_{2}}\| 	\\
& \lesssim \|\phi\|_{\mathcal{B}^{-\frac{1}{4},\frac{1}{2};\infty}_{\pm}}\sum_{\mathbf{N,L}}N_{0}^{1/4}N_{1}^{-1/2}\left(\frac{L_{0}}{N_{1}}\right)^{1/2}(N_{1}L_{1})^{1/2}(N_{1}L_{2})^{1/4}\|P_{K_{N_{1},L_{1}}^{\pm_{1}}}\psi_{1,\pm_{1}}\|\|P_{K_{N_{2},L_{2}}^{\pm_{2}}}\psi_{2,\pm_{2}}\| \\
&= \|\phi\|_{\mathcal{B}^{-\frac{1}{4},\frac{1}{2};\infty}_{\pm}}\sum_{\mathbf{N,L}}N_{0}^{1/4}\left(\frac{L_{0}}{N_{1}}\right)^{1/2}L_{1}^{1/2}(N_{1}L_{2})^{1/4}\|P_{K_{N_{1},L_{1}}^{\pm_{1}}}\psi_{1,\pm_{1}}\|\|P_{K_{N_{2},L_{2}}^{\pm_{2}}}\psi_{2,\pm_{2}}\| \\
& \lesssim \|\psi_{1,\pm_{1}}\|_{\mathcal{B}^{\frac{1}{4},\frac{1}{2};1}_{\pm_{1}}}\|\psi_{2,\pm_{2}}\|_{\mathcal{B}^{\frac{1}{4},\frac{1}{2};1}_{\pm_{2}}}\|\phi\|_{\mathcal{B}^{-\frac{1}{4},\frac{1}{2};\infty}_{\pm}}.
\end{align*}

Finally, we are left to treat the case $L_2\sim L_0$. This is very straightforward.
If $N_{0}\ll N_{1}\sim N_{2}$, then we use the bilinear estimate \eqref{bi-selb-1} to get
\begin{align*}
|\mathbf J^1| &\lesssim \|\phi\|_{\mathcal{B}^{-\frac{1}{4},\frac{1}{2};\infty}_{\pm}}\sum_{\mathbf{N,L}}N_{0}^{1/4}L_{0}^{-1/2}\bigg(\frac{L_{0}}{N_{1}}\bigg)^{1/2}(N_{0}L_{1})^{1/2}(N_{1}L_{2})^{1/4}\|P_{K_{N_{1},L_{1}}^{\pm_{1}}}\psi_{1,\pm_{1}}\|\|P_{K_{N_{2},L_{2}}^{\pm_{2}}}\psi_{2,\pm_{2}}\| 	\\
& = \|\phi\|_{\mathcal{B}^{-\frac{1}{4},\frac{1}{2};\infty}_{\pm}}\sum_{\mathbf{N,L}}\bigg(\frac{N_{0}}{N_{1}}\bigg)^{1/2}N_{0}^{1/4}L_{1}^{1/2}N_{1}^{1/4}L_{2}^{1/4}\|P_{K_{N_{1},L_{1}}^{\pm_{1}}}\psi_{1,\pm_{1}}\|\|P_{K_{N_{2},L_{2}}^{\pm_{2}}}\psi_{2,\pm_{2}}\|\\
& \lesssim  \|\psi_{1,\pm_{1}}\|_{\mathcal{B}^{\frac{1}{4},\frac{1}{2};1}_{\pm_{1}}}\|\psi_{2,\pm_{2}}\|_{\mathcal{B}^{\frac{1}{4},\frac{1}{2};1}_{\pm_{2}}}\|\phi\|_{\mathcal{B}^{-\frac{1}{4},\frac{1}{2};\infty}_{\pm}} \;\;\left(\because \sum_{L_{0};L_{0}\sim L_{2}}\sum_{N_{0};N_{0}\ll N_{1}}\bigg(\frac{N_{0}}{N_{1}}\bigg)^{1/2}\lesssim 1\right).
\end{align*}

Now we consider the case $N_{1}\lesssim N_{0}\sim N_{2}$ and estimate
\begin{align*}
|\mathbf J^1| &\lesssim \|\phi\|_{\mathcal{B}^{-\frac{1}{4},\frac{1}{2};\infty}_{\pm}}\sum_{\mathbf{N,L}}N_{0}^{1/4}L_{0}^{1/2}\bigg(\frac{L_{0}}{N_{1}}\bigg)^{1/2}(N_{1}L_{1})^{1/2}(N_{1}L_{2})^{1/4}\|P_{K_{N_{1},L_{1}}^{\pm_{1}}}\psi_{1,\pm_{1}}\|\|P_{K_{N_{2},L_{2}}^{\pm_{2}}}\psi_{2,\pm_{2}}\|\\
& = \|\phi\|_{\mathcal{B}^{-\frac{1}{4},\frac{1}{2};\infty}_{\pm}}\sum_{\mathbf{N,L}}N_{0}^{1/4}L_{1}^{1/2}N_{1}^{1/4}L_{2}^{1/4}\|P_{K_{N_{1},L_{1}}^{\pm_{1}}}\psi_{1,\pm_{1}}\|\|P_{K_{N_{2},L_{2}}^{\pm_{2}}}\psi_{2,\pm_{2}}\| \\
& \lesssim  \|\psi_{1,\pm_{1}}\|_{\mathcal{B}^{\frac{1}{4},\frac{1}{2};1}_{\pm_{1}}}\|\psi_{2,\pm_{2}}\|_{\mathcal{B}^{\frac{1}{4},\frac{1}{2},1}_{\pm_{2}}}\|\phi\|_{\mathcal{B}^{-\frac{1}{4},\frac{1}{2};\infty}_{\pm}}\;\;\left(\because \sum_{L_{0};L_{0}\sim L_{2}}\sum_{N_{0};N_{0}\sim N_{2}}N_{0}^{1/4}\lesssim N_{2}^{1/4}\right).	
\end{align*}

\subsection{Proof of \eqref{nl-anu-2}}
We write
\begin{eqnarray*}
\|R^{\mu}_{\pm}\epsilon_{\mu\nu\lambda}\psi_{1,\pm_{1}}^{\dagger}R_{\pm_{2}}^{\lambda}\psi_{2,\pm_{2}}\|_{\mathcal{B}^{\frac{1}{4},-\frac{1}{2};1}_{\pm}} & = & \sup_{\|\phi\|_{\mathcal{B}^{-\frac{1}{4},\frac{1}{2};\infty}_{\pm}}=1} \left| \int \epsilon_{\mu\nu\lambda}\psi_{1,\pm_{1}}^{\dagger}R^{\lambda}_{\pm_{2}}\psi_{2,\pm_{2}}\overline{R^{\mu}_{\pm}\phi}\,dtdx \right|	
\end{eqnarray*}
and set
$$
\mathbf J^2 :=  \int \epsilon_{\mu\nu\lambda}\psi_{1,\pm_{1}}^{\dagger}R^{\lambda}_{\pm_{2}}\psi_{2,\pm_{2}}\overline{R^{\mu}_{\pm}\phi}\,dtdx.
$$
Then the dyadic decomposition gives $|\mathbf J^2| \le \sum_{\mathbf{N,L}}|\mathbf J^2_{\mathbf{N,L}}|$, where
\begin{align*}
|\mathbf J^2_{\mathbf{N,L}}| &=  \left|\int \epsilon_{\mu\nu\lambda} P_{K_{N_{1},L_{1}}^{\pm_{1}}}\psi_{1,\pm_{1}}^{\dagger}R^{\lambda}_{\pm_{2}}P_{K_{N_{2},L_{2}}^{\pm_{2}}}\psi_{2,\pm_{2}}\overline{R^{\mu}_{\pm}P_{K_{N_{0},L_{0}}^{\pm}}\phi}\,dtdx\right|\\
&=  \left|\int \mathcal F( P_{K_{N_{1},L_{1}}^{\pm_{1}}}\psi_{1,\pm_{1}})^{\dagger}(-X_1)\epsilon_{\mu\nu\lambda} \mathcal F( R^{\lambda}_{\pm_{2}}P_{K_{N_{2},L_{2}}^{\pm_{2}}}\psi_{2,\pm_{2}}\overline{R^{\mu}_{\pm}P_{K_{N_{0},L_{0}}^{\pm}}\phi})(X_1)\,dX_1\right|.
\end{align*}
Also we can represent $\mathbf J^2_{\mathbf{N,L}}$ by
$$
|\mathbf J^2_{\mathbf{N,L}}| = \left|\int\mathcal F(P_{K_{N_1,L_1}^{\pm_1}}\psi_{1,\pm_1}^\dagger R_{\pm_2}^\lambda P_{K_{N_2,L_2}^{\pm_2}}\psi_{2,\pm_2})(X_1)\overline{\mathcal F (R_\pm^\mu P_{K_{N_0,L_0}^\pm}\phi)}(-X_1)dX_1\right|.
$$

\begin{rem}\label{hll-mod}
Now we give same argument by Remark \ref{llh-mod}. Indeed, if $L_{\max}^{02}\ll L_1$ and $N_1\ll N_0\sim N_2$, we see that the support condition excludes $\pm\neq\pm_2$ in view of first representation of $\mathbf J^2_{\mathbf{N,L}}$. Also, from second representation, we excludes the case $\pm=\pm_2$. This implies that llh-hhl case does not appear in the summation of $\mathbf J^2_{\mathbf{N,L}}$. On the other hand, for $N_0\ll N_1\sim N_2$, first representation gives the exclusion of the case $\pm_1=\pm_2$ and second gives $\pm_1\neq\pm_2$. We can argue similarly in the case $N_2\ll N_0\sim N_1$ and hence the case llh-hlh does not appear. By the same reason the case $L_{\max}^{01}\ll L_2\ll N_2\ll N_0\sim N_1$ can be also excluded from the summation. In fact, from the first representation of $\mathbf J^2_{\mathbf{N,L}}$, we exclude the case $\pm=\pm_1$ and the case $\pm\neq\pm_2$ from the second representation.
\end{rem}

Since $\overline{R^{\mu}_{\pm}P_{K_{N_{0},L_{0}}^{\pm}}\phi} = R_{\pm}^\mu(-\overline{P_{K_{N_{0},L_{0}}^{\pm}}\phi})$, we write
$$
|\mathbf J^2| = \sum_{\mathbf{N,L}} \int \mathcal F( P_{K_{N_{1},L_{1}}^{\pm_{1}}}\psi_{1,\pm_{1}})^{\dagger}\epsilon_{\mu\nu\lambda} \mathcal F[ Q^{\lambda\mu}_{\pm_{2},\pm}(P_{K_{N_{2},L_{2}}^{\pm_{2}}}\psi_{2,\pm_{2}}, -\overline{P_{K_{N_{0},L_{0}}^{\pm}}\phi})]\,dX_1.
$$
By the same reason as in the previous section we may assume that $\mathcal F P_{K_{N_{0},L_{0}}^{\pm}}\phi$, $\mathcal F P_{K_{N_{1},L_{1}}^{\pm_{1}}}\psi_{1,\pm_{1}}$, and $\mathcal F P_{K_{N_{2},L_{2}}^{\pm_{2}}}\psi_{2,\pm_{2}}$ are nonnegative real-valued functions.
Thus by using \eqref{bil-1} and Lemma \ref{bi-int}, we get
\begin{eqnarray*}
|\mathbf J^2| & \lesssim & \sum_{\mathbf{N,L}}\int \mathcal F P_{K_{N_{1},L_{1}}^{\pm_{1}}}\psi_{1,\pm_{1}}(X_{1})\int_{X_{1}=X_{0}-X_{2}}\theta_{02}\;\mathcal F P_{K_{N_{0},L_{0}}^{\pm}}\phi(X_{0})
\mathcal F P_{K_{N_{2},L_{2}}^{\pm_{2}}}\psi_{2,\pm_{2}}(X_{2})\,dX_{0}dX_{2}dX_{1}\\
& \lesssim & \sum_{\mathbf{N,L}}\|P_{K_{N_{1},L_{1}}^{\pm_{1}}}\psi_{1,\pm_{1}}\|\bigg(\frac{L_{\max}^{012}}{N_{\min}^{02}}\bigg)^{1/2}
\|P_{K_{N_{1},L_{1}}^{\pm_{1}}}(P_{K_{N_{0},L_{0}}^{\pm}}\phi P_{K_{N_{2},L_{2}}^{\pm_{2}}}\psi_{2,\pm_{2}})\|.	
\end{eqnarray*}

\subsubsection{Case 1: $L_{1}\le L_{0}\le L_{2}$}\label{l102}\hfill

 If $N_{1}\ll N_{0}\sim N_{2}$, then
by \eqref{bi-selb-2} we get
\begin{align*}
	|\mathbf J^2| & \lesssim \|\phi\|_{\mathcal{B}^{-\frac{1}{4},\frac{1}{2};\infty}_{\pm}}\sum_{\mathbf{N,L}}N_{0}^{1/4}L_{0}^{-1/2}\left(\frac{L_2}{N_0}\right)^{1/2}(N_1L_1)^{1/2}(N_0L_0)^{1/4}\|P_{K_{N_{1},L_{1}}^{\pm_{1}}}\psi_{1,\pm_{1}}\|\;\|P_{K_{N_{2},L_{2}}^{\pm_{2}}}\psi_{2,\pm_{2}}\| \\
	&\lesssim  \|\phi\|_{\mathcal{B}^{-\frac{1}{4},\frac{1}{2};\infty}_{\pm}}\sum_{\mathbf{N,L}}N_1^{1/2}L_0^{-1/4}L_1^{1/2}L_2^{1/2}\|P_{K_{N_{1},L_{1}}^{\pm_{1}}}\psi_{1,\pm_{1}}\|\;\|P_{K_{N_{2},L_{2}}^{\pm_{2}}}\psi_{2,\pm_{2}}\| \\
	& \le \|\psi_{1,\pm_{1}}\|_{\mathcal{B}^{\frac{1}{4},\frac{1}{2};1}_{\pm_{1}}}\|\psi_{2,\pm_{2}}\|_{\mathcal{B}^{\frac{1}{4},\frac{1}{2};1}_{\pm_{2}}}\|\phi\|_{\mathcal{B}^{-\frac{1}{4},\frac{1}{2};\infty}_{\pm}}\;\;\left(\because \sum_{L_{0}\ge1}\sum_{N_{0};N_{0}\gg N_{1}}N_{0}^{-1/4}L_{0}^{-1/4}\lesssim N_{1}^{-1/4} \right).
\end{align*}

If $N_{0}\lesssim N_{1}\sim N_{2}$, then \eqref{bi-selb-2} gives
\begin{align*}
|\mathbf J^2|& \lesssim \|\phi\|_{\mathcal{B}^{-\frac{1}{4},\frac{1}{2};\infty}_{\pm}}\sum_{\mathbf{N,L}}N_{0}^{1/4}L_{0}^{-1/2}\left(\frac{L_2}{N_0}\right)^{1/2}(N_0L_1)^{1/2}(N_0L_0)^{1/4}\|P_{K_{N_{1},L_{1}}^{\pm_{1}}}\psi_{1,\pm_{1}}\|\;\|P_{K_{N_{2},L_{2}}^{\pm_{2}}}\psi_{2,\pm_{2}}\| \\
&\lesssim \|\phi\|_{\mathcal{B}^{-\frac{1}{4},\frac{1}{2};\infty}_{\pm}}\sum_{\mathbf{N,L}}N_0^{1/2}L_0^{-1/4}L_1^{1/2}L_2^{1/2}\|P_{K_{N_{1},L_{1}}^{\pm_{1}}}\psi_{1,\pm_{1}}\|\;\|P_{K_{N_{2},L_{2}}^{\pm_{2}}}\psi_{2,\pm_{2}}\| \\
& \le \|\psi_{1,\pm_{1}}\|_{\mathcal{B}^{\frac{1}{4},\frac{1}{2};1}_{\pm_{1}}}\|\psi_{2,\pm_{2}}\|_{\mathcal{B}^{\frac{1}{4},\frac{1}{2};1}_{\pm_{2}}}\|\phi\|_{\mathcal{B}^{-\frac{1}{4},\frac{1}{2};\infty}_{\pm}}\;\;\left(\because \sum_{L_{0}\ge1}\sum_{N_{0};N_{0}\lesssim N_{1}}N_{0}^{1/2}L_{0}^{-1/4}\lesssim N_{1}^{1/2}\lesssim N_{1}^{1/4}N_{2}^{1/4}\right).
\end{align*}

If $N_{2}\lesssim N_{0}\sim N_{1}$ and $L_2\gtrsim N_2$, then we see that $L_2\ll N_1\sim N_0$ is excluded by the same way as Remark \ref{hll-mod}. Then we only need to treat the case $L_2\sim N_0\sim N_1$. Note that in the case $L_2\gtrsim N_2$, we have $\theta_{02}\sim1$, and hence instead of using $\dfrac{L_2}{N_2}$, we use $\dfrac{L_2}{N_0}\sim1$ and \eqref{bi-selb-2} with $j=2$ to get
\begin{align*}
|\mathbf J^2| &\lesssim \|\phi\|_{\mathcal{B}^{-\frac{1}{4},\frac{1}{2};\infty}_{\pm}}\sum_{\mathbf{N,L}}N_0^{1/4}L_0^{-1/2}\left(\frac{L_2}{N_0}\right)^{1/2}(N_2L_1)^{1/2}(N_0L_0)^{1/4} \|P_{K_{N_1,L_1}^{\pm_1}}\psi_{1,\pm_1}\|\|P_{K_{N_2,L_2}^{\pm_2}}\psi_{2,\pm_2}\| \\
&\lesssim \|\phi\|_{\mathcal{B}^{-\frac{1}{4},\frac{1}{2};\infty}_{\pm}}\sum_{\mathbf{N,L}}N_2^{1/2}L_0^{-1/4}L_1^\frac12 L_2^\frac12\|P_{K_{N_1,L_1}^{\pm_1}}\psi_{1,\pm_1}\|\|P_{K_{N_2,L_2}^{\pm_2}}\psi_{2,\pm_2}\|\\
& \lesssim \|\psi_{1,\pm_{1}}\|_{\mathcal{B}^{\frac{1}{4},\frac{1}{2};1}_{\pm_{1}}}\|\psi_{2,\pm_{2}}\|_{\mathcal{B}^{\frac{1}{4},\frac{1}{2};1}_{\pm_{2}}}\|\phi\|_{\mathcal{B}^{-\frac{1}{4},\frac{1}{2};\infty}_{\pm}}.
\end{align*}

Now let us consider the case $L_{\max}^{01}\ll L_2\ll N_2$. By the exclusion as stated in Remark \ref{hll-mod} we have only to consider the case $L_{\max}^{01}\ll L_2\ll N_2\sim N_0\sim N_1$. To do this, we follow the argument in Case 3 in the Section 6.1. First, since $N_2\sim N_1\sim N_0$, by the law of sine, we deduce that $\theta_{02}\lesssim \dfrac{N_0}{N_2}\theta_{01}\sim\theta_{01}$ and then we estimate $\mathbf J^2_{\mathbf{N,L}}$ as
\begin{align*}
|\mathbf J^2_{\mathbf{N,L}}| &\lesssim \int\!\int_{X_1=X_2-X_0} \theta_{02}\mathcal F(P_{K_{N_1,L_1}^{\pm_1}}\psi_{1,\pm_1})^\dagger(X_1)F(P_{K_{N_0,L_0}^{\pm}}\phi)(X_0)F(P_{K_{N_2,L_1}^{\pm_2}}\psi_{2,\pm_2})(X_2)dX_2dX_0dX_1\\
&\sim \int\!\int_{X_2=X_1+X_0} \theta_{01}\mathcal F(P_{K_{N_1,L_1}^{\pm_1}}\psi_{1,\pm_1})^\dagger(X_1)F(P_{K_{N_0,L_0}^{\pm}}\phi)(X_0)F(P_{K_{N_2,L_1}^{\pm_2}}\psi_{2,\pm_2})(X_2)dX_1dX_0dX_2.
\end{align*}

Using \eqref{ang-whi2}, we write
$$
|\mathbf J^2| \lesssim \sum_{\mathbf{N,L}}\sum_{\omega_1,\omega_0}\int \mathcal F [P_{K_{N_1,L_1}^{\pm_1}}\psi_{1,\pm_1}^{2\theta_{01},\omega_1}\overline{R_\pm^\mu P_{K_{N_0,L_0}^\pm}\phi^{2\theta_{01},\omega_0}}]\mathcal F [R_{\pm_2}^\lambda P_{K_{N_2,L_2}^{\pm_2}}\psi_{2,\pm_2}]dX_2,
$$
where the spatial Fourier support of $P_{K_{N_1,L_1}^{\pm_1}}\psi_{1,\pm_1}^{2\theta_{01},\omega_1}$ is contained in a strip of radius compatible to $N_{\max}^{01}\theta_{01}$ around $\mathbb R\omega_1$.

Then by Cauchy-Schwarz inequality and Theorem \ref{null-form-selb} with $r\sim N_{\max}^{01}\theta_{01}$, we have
\begin{align*}
|\mathbf J^2| &\lesssim \sum_{\mathbf{N,L}}\sum_{\omega_1,\omega_0}\|B_{\theta_{01}}(P_{K_{N_1,L_1}^{\pm_1}}\psi_{1,\pm_1}^{2\theta_{01},\omega_1},R_\pm^\mu P_{K_{N_0,L_0}^\pm}\phi^{2\theta_{01},\omega_0})\|\|P_{K_{N_2,L_2}^{\pm_2}}\psi_{2,\pm_2}\| \\
&\lesssim \|\phi\|_{\mathcal{B}^{-\frac{1}{4},\frac{1}{2};\infty}_{\pm}}\sum_{\mathbf{N,L}}N_0^{1/4}L_0^{-1/2}\left(N_{\max}^{01}\left(\frac{L_2}{N_{\min}^{01}}\right)^{1/2}L_0L_1\right)^{1/2}\|P_{K_{N_1,L_1}^{\pm_1}}\psi_{1,\pm_1}\|\|P_{K_{N_2,L_2}^{\pm_2}}\psi_{2,\pm_2}\| \\
&\lesssim \|\phi\|_{\mathcal{B}^{-\frac{1}{4},\frac{1}{2};\infty}_{\pm}}\sum_{\mathbf{N,L}}N_0^{1/2}L_1^{1/2}L_2^{1/4}\|P_{K_{N_1,L_1}^{\pm_1}}\psi_{1,\pm_1}\|\|P_{K_{N_2,L_2}^{\pm_2}}\psi_{2,\pm_2}\|\\
& \lesssim \|\psi_{1,\pm_{1}}\|_{\mathcal{B}^{\frac{1}{4},\frac{1}{2};1}_{\pm_{1}}}\|\psi_{2,\pm_{2}}\|_{\mathcal{B}^{\frac{1}{4},\frac{1}{2};1}_{\pm_{2}}}\|\phi\|_{\mathcal{B}^{-\frac{1}{4},\frac{1}{2};\infty}_{\pm}},
\end{align*}
where the summation by $\omega_1,\omega_0\in\Omega(\theta_{01})$ is taken over $\angle(\omega_1,\omega_0)\le4\theta_{01}$.

\subsubsection{Case 2: $L_{0} \le L_{1} \le L_{2}$}\hfill

For $N_{1}\ll N_{0}\sim N_{2}$, we have
\begin{align*}
|\mathbf J^2| &\lesssim \|\phi\|_{\mathcal{B}^{-\frac{1}{4},\frac{1}{2};\infty}_{\pm}}\sum_{\mathbf{N,L}}N_{0}^{1/4}L_{0}^{-1/2}\left(\frac{L_2}{N_0}\right)^{1/2}(N_1L_0)^{1/2}(N_1L_1)^{1/4}\|P_{K_{N_{1},L_{1}}^{\pm_{1}}}\psi_{1,\pm_{1}}\|\;\|P_{K_{N_{2},L_{2}}^{\pm_{2}}}\psi_{2,\pm_{2}}\| \\
&\lesssim \|\phi\|_{\mathcal{B}^{-\frac{1}{4},\frac{1}{2};\infty}_{\pm}}\sum_{\mathbf{N,L}}N_0^{1/4}N_1^{1/4}L_1^{1/4}L_2^{1/2}\|P_{K_{N_{1},L_{1}}^{\pm_{1}}}\psi_{1,\pm_{1}}\|\;\|P_{K_{N_{2},L_{2}}^{\pm_{2}}}\psi_{2,\pm_{2}}\| \\
 & \lesssim  \|\psi_{1,\pm_{1}}\|_{\mathcal{B}^{\frac{1}{4},\frac{1}{2};1}_{\pm_{1}}}\|\psi_{2,\pm_{2}}\|_{\mathcal{B}^{\frac{1}{4},\frac{1}{2};1}_{\pm_{2}}}\|\phi\|_{\mathcal{B}^{-\frac{1}{4},\frac{1}{2};\infty}_{\pm}}\left(\because \sum_{L_{0};L_{0}\le L_{1}}\sum_{N_{0};N_{0}\gg N_{1}}N_{0}^{-1/4}\lesssim N_{1}^{-1/4}L_{1}^{\epsilon}\right).
\end{align*}

If $N_{0} \lesssim N_{1}\sim N_{2}$, then
\begin{align*}
|\mathbf J^2| &\lesssim \|\phi\|_{\mathcal{B}^{-\frac{1}{4},\frac{1}{2};\infty}_{\pm}}\sum_{\mathbf{N,L}}N_{0}^{1/4}L_{0}^{-1/2}\left(\frac{L_2}{N_0}\right)^{1/2}(N_0L_0)^{1/2}(N_0L_1)^{1/4}\|P_{K_{N_{1},L_{1}}^{\pm_{1}}}\psi_{1,\pm_{1}}\|\|P_{K_{N_{2},L_{2}}^{\pm_{2}}}\psi_{2,\pm_{2}}\|\\
&\lesssim \|\phi\|_{\mathcal{B}^{-\frac{1}{4},\frac{1}{2};\infty}_{\pm}}\sum_{\mathbf{N,L}}N_0^{1/2}L_1^{1/4}L_2^{1/2}\|P_{K_{N_{1},L_{1}}^{\pm_{1}}}\psi_{1,\pm_{1}}\|\|P_{K_{N_{2},L_{2}}^{\pm_{2}}}\psi_{2,\pm_{2}}\|\\
& \lesssim  \|\psi_{1,\pm_{1}}\|_{\mathcal{B}^{\frac{1}{4},\frac{1}{2};1}_{\pm_{1}}}\|\psi_{2,\pm_{2}}\|_{\mathcal{B}^{\frac{1}{4},\frac{1}{2};1}_{\pm_{2}}}\|\phi\|_{\mathcal{B}^{-\frac{1}{4},\frac{1}{2};\infty}_{\pm}}	\left(\because \sum_{L_{0};L_{0}\le L_{1}}\sum_{N_{0};N_{0}\lesssim N_{1}}N_{0}^{1/2}\lesssim N_{1}^{1/2}L_{1}^{\epsilon}\lesssim N_{1}^{1/4}N_{2}^{1/4}L_{1}^{\epsilon}\right).
\end{align*}

If $N_{2} \lesssim N_{0}\sim N_{1}$, we can exclude the case $L_{\max}^{01}\ll L_2\ll N_2\ll N_0\sim N_1$ by the argument of Remark \ref{hll-mod}. The remaining cases can be treated similarly to the last case of Section \ref{l102} by changing the role of $N_1, N_2$. We omit the details.

\subsubsection{Case 3: $L_{0}\le L_{2}\le L_{1}$}\hfill

%Similarly, we first consider the case $L_2\ll L_1$. As the proof of \eqref{nl-anu-1}, we see that the supports are disjoint if $L_1\ll N_{\min}^{012}$ with suitable choices of the signs $\pm_j'$. Also, for the case $L_1\lesssim N_{\min}^{012}$, we get the required estimtes directly from the inequality $L_1\lesssim N_{\min}^{012}$. Thus we are left to treat the case $L_2\sim L_1$ which is straightforward.

For the low-low-high modulation, we remains to consider the case $N_0\sim N_1\sim N_2$.  The angle Whitney decomposition gives us
$$
|\mathbf J^2| \lesssim \sum_{\mathbf{N,L}}\sum_{\omega_2,\omega_0}\int \mathcal F [(P_{K_{N_1,L_1}^{\pm_1}}\psi_{1,\pm_1})^\dagger](-X_1) \theta_{02}\mathcal F[ P_{K_{N_2,L_2}^{\pm_2}}\psi_{2,\pm_2}^{2\theta_{02},\omega_2}\overline{ P_{K_{N_0,L_0}^{\pm}}\phi^{2\theta_{02},\omega_0}}](X_1)dX_1.
$$
Here, the spatial Fourier support of $P_{K_{N_2,L_2}^{\pm_2}}\psi_{2,\pm_2}^{2\theta_{02},\omega_2}$ is contained in a strip of radius compatible to $N_{\max}^{02}\theta_{02}$ around $\mathbb R\omega_2$.

Then by Theorem \ref{null-form-selb} with $r\sim N_{\max}^{02}\theta_{02}$,
\begin{align*}
|\mathbf J^2| &\lesssim \sum_{\mathbf{N,L}}\sum_{\omega_2,\omega_0}\|P_{K_{N_1,L_1}^{\pm_1}}\psi_{1,\pm_1}\|\|P_{K_{N_1,L_1}^{\pm_1}}B_{\theta_{02}}(P_{K_{N_2,L_2}^{\pm_2}}\psi_{2,\pm_2}^{2\theta_{02},\omega_2},P_{K_{N_0,L_0}^{\pm}}\phi^{2\theta_{02},\omega_0})\| \\
&\lesssim \sum_{\mathbf{N,L}}\sum_{\omega_2,\omega_0}\|P_{K_{N_1,L_1}^{\pm_1}}\psi_{1,\pm_1}\|\left(N_{\max}^{02}\left(\frac{L_1}{N_{\min}^{02}}\right)^{1/2}L_0L_2\right)^{1/2}\|P_{K_{N_2,L_2}^{\pm_2}}\psi_{2,\pm_2}^{2\theta_{02},\omega_2}\|\|P_{K_{N_0,L_0}^{\pm}}\phi^{2\theta_{02},\omega_0}\| \\
& \lesssim \|\phi\|_{\mathcal{B}^{-\frac{1}{4},\frac{1}{2};\infty}_{\pm}}\sum_{\mathbf{N,L}}N_0^{1/2}L_1^{1/4}L_2^{1/2}\|P_{K_{N_1,L_1}^{\pm_1}}\psi_{1,\pm_1}\|\|P_{K_{N_2,L_2}^{\pm_2}}\psi_{2,\pm_2}\| \\
&\lesssim \|\psi_{1,\pm_{1}}\|_{\mathcal{B}^{\frac{1}{4},\frac{1}{2};1}_{\pm_{1}}}\|\psi_{2,\pm_{2}}\|_{\mathcal{B}^{\frac{1}{4},\frac{1}{2};1}_{\pm_{2}}}\|\phi\|_{\mathcal{B}^{-\frac{1}{4},\frac{1}{2};\infty}_{\pm}}.
\end{align*}

Now we consider the case $L_1\gtrsim N_{\min}^{012}$. If $N_1\ll N_2\sim N_0$, then
\begin{align*}
|\mathbf J^2| &\lesssim \|\phi\|_{\mathcal{B}^{-\frac{1}{4},\frac{1}{2};\infty}_{\pm}}\sum_{\mathbf{N,L}}N_0^{1/4}L_0^{-1/2}\|P_{K_{N_1,L_1}^{\pm_1}}\psi_{1,\pm_1}\|\left(\frac{L_1}{N_0}\right)^{1/2}(N_1L_0)^{1/2}(N_0L_2)^{1/4}\|P_{K_{N_2,L_2}^{\pm_2}}\psi_{2,\pm_2}\| \\
& \lesssim \|\phi\|_{\mathcal{B}^{-\frac{1}{4},\frac{1}{2};\infty}_{\pm}}\sum_{\mathbf{N,L}}N_1^{1/2}L_1^{1/2}L_2^{1/4}\|P_{K_{N_1,L_1}^{\pm_1}}\psi_{1,\pm_1}\|\|P_{K_{N_2,L_2}^{\pm_2}}\psi_{2,\pm_2}\| \\
&\lesssim \|\psi_{1,\pm_{1}}\|_{\mathcal{B}^{\frac{1}{4},\frac{1}{2};1}_{\pm_{1}}}\|\psi_{2,\pm_{2}}\|_{\mathcal{B}^{\frac{1}{4},\frac{1}{2};1}_{\pm_{2}}}\|\phi\|_{\mathcal{B}^{-\frac{1}{4},\frac{1}{2};\infty}_{\pm}}.
\end{align*}

For $N_0\ll N_1\sim N_2$, we have
\begin{align*}
|\mathbf J_2| &\lesssim \|\phi\|_{\mathcal{B}^{-\frac{1}{4},\frac{1}{2};\infty}_{\pm}}\sum_{\mathbf{N,L}}N_0^{1/4}L_0^{-1/2}\|P_{K_{N_1,L_1}^{\pm_1}}\psi_{1,\pm_1}\|\left(\frac{L_1}{N_0}\right)^{1/2}(N_0L_0)^{1/2}(N_0L_2)^{1/4}\|P_{K_{N_2,L_2}^{\pm_2}}\psi_{2,\pm_2}\| \\
&\lesssim \|\phi\|_{\mathcal{B}^{-\frac{1}{4},\frac{1}{2};\infty}_{\pm}}\sum_{\mathbf{N,L}}N_0^{1/2}L_1^{1/2}L_2^{1/4}\|P_{K_{N_1,L_1}^{\pm_1}}\psi_{1,\pm_1}\|\|P_{K_{N_2,L_2}^{\pm_2}}\psi_{2,\pm_2}\| \\
&\lesssim \|\psi_{1,\pm_{1}}\|_{\mathcal{B}^{\frac{1}{4},\frac{1}{2};1}_{\pm_{1}}}\|\psi_{2,\pm_{2}}\|_{\mathcal{B}^{\frac{1}{4},\frac{1}{2};1}_{\pm_{2}}}\|\phi\|_{\mathcal{B}^{-\frac{1}{4},\frac{1}{2};\infty}_{\pm}}.
\end{align*}

If $N_2\ll N_1\sim N_0$, then
\begin{align*}
|\mathbf J_2| &\lesssim \|\phi\|_{\mathcal{B}^{-\frac{1}{4},\frac{1}{2};\infty}_{\pm}}\sum_{\mathbf{N,L}}N_0^{1/4}L_0^{-1/2}\|P_{K_{N_1,L_1}^{\pm_1}}\psi_{1,\pm_1}\|\left(\frac{L_1}{N_2}\right)^{1/2}(N_2L_0)^{1/2}(N_2L_2)^{1/4}\|P_{K_{N_2,L_2}^{\pm_2}}\psi_{2,\pm_2}\| \\
&\lesssim \|\phi\|_{\mathcal{B}^{-\frac{1}{4},\frac{1}{2};\infty}_{\pm}}\sum_{\mathbf{N,L}}N_0^{1/4}N_2^{1/4}L_1^{1/2}L_2^{1/4}\|P_{K_{N_1,L_1}^{\pm_1}}\psi_{1,\pm_1}\|\|P_{K_{N_2,L_2}^{\pm_2}}\psi_{2,\pm_2}\|\\
&\lesssim \|\psi_{1,\pm_{1}}\|_{\mathcal{B}^{\frac{1}{4},\frac{1}{2};1}_{\pm_{1}}}\|\psi_{2,\pm_{2}}\|_{\mathcal{B}^{\frac{1}{4},\frac{1}{2};1}_{\pm_{2}}}\|\phi\|_{\mathcal{B}^{-\frac{1}{4},\frac{1}{2};\infty}_{\pm}}.
\end{align*}

Finally it remains to deal with the case $L_2\sim L_1$. Here, we only consider the case $N_{1}\ll N_{0}\sim N_{2}$. The other cases can be treated similarly.
\begin{align*}
|\mathbf J^2| &\lesssim \|\phi\|_{\mathcal{B}^{-\frac{1}{4},\frac{1}{2};\infty}_{\pm}}\sum_{\mathbf{N,L}}N_{0}^{1/4}L_{0}^{-1/2}\bigg(\frac{L_{1}}{N_{0}}\bigg)^{1/2}(N_{1}L_{0})^{1/2}(N_{0}L_{2})^{1/4}\|P_{K_{N_{1},L_{1}}^{\pm_{1}}}\psi_{1,\pm_{1}}\|\;\|P_{K_{N_{2},L_{2}}^{\pm_{2}}}\psi_{2,\pm_{2}}\| \\
& =  \|\phi\|_{\mathcal{B}^{-\frac{1}{4},\frac{1}{2};\infty}_{\pm}}\sum_{\mathbf{N,L}}L_{1}^{1/2}N_{1}^{1/2}L_{2}^{1/4}\|P_{K_{N_{1},L_{1}}^{\pm_{1}}}\psi_{1,\pm_{1}}\|\;\|P_{K_{N_{2},L_{2}}^{\pm_{2}}}\psi_{2,\pm_{2}}\|\\
& \lesssim  \|\phi\|_{\mathcal{B}^{-\frac{1}{4},\frac{1}{2};\infty}_{\pm}}\sum_{\mathbf{N,L}}N_{1}^{1/4}L_{1}^{1/4}N_{0}^{1/4}L_{2}^{1/4}\|P_{K_{N_{1},L_{1}}^{\pm_{1}}}\psi_{1,\pm_{1}}\|\;\|P_{K_{N_{2},L_{2}}^{\pm_{2}}}\psi_{2,\pm_{2}}\|\\
& \lesssim \|\psi_{1,\pm_{1}}\|_{\mathcal{B}^{\frac{1}{4},\frac{1}{2};1}_{\pm_{1}}}\|\psi_{2,\pm_{2}}\|_{\mathcal{B}^{\frac{1}{4},\frac{1}{2};1}_{\pm_{2}}}\|\phi\|_{\mathcal{B}^{-\frac{1}{4},\frac{1}{2};\infty}_{\pm}} \left(\because \sum_{L_{1};L_{1}\le L_{2}}\sum_{N_{1};N_{1}\sim N_{2}}N_{1}^{1/4}\lesssim N_{2}^{1/4}L_{2}^{\epsilon}\right).
\end{align*}

This completes the proof of \eqref{psipsi}.

\section{Estimates of $\mathcal M$: Proof of \eqref{psia}}

We need to show that
\begin{align}
\|\Pi_{\pm_{2}}(A_{\mu}\Pi_{\mp_{1}}\alpha^{\mu}\Pi_{\pm_{1}}\psi_{1,\pm_{1}})\|_{\mathcal{B}^{\frac{1}{4},-\frac{1}{2};1}_{\pm_{2}}}	& \lesssim  \|A_{\mu,\pm_{0}}\|_{\mathcal{B}^{\frac{1}{4},\frac{1}{2};1}_{\pm_{0}}}\|\psi_{1,\pm_{1}}\|_{\mathcal{B}^{\frac{1}{4},\frac{1}{2};1}_{\pm_{1}}},\label{psia1} \\
\|\Pi_{\pm_{2}}(A_{\mu}R^{\mu}_{\pm}\psi_{\pm})\|_{\mathcal{B}^{\frac{1}{4},-\frac{1}{2};1}_{\pm_{2}}} & \lesssim  \|A_{\mu,\pm}\|_{\mathcal{B}^{\frac{1}{4},\frac{1}{2};1}_{\pm}}\|\psi_{\pm}\|_{\mathcal{B}^{\frac{1}{4},\frac{1}{2};1}_{\pm}}.\label{psia2}
\end{align}
The estimate \eqref{psia1} can be done by duality, dyadic decomposition, and bilinear estimates without any essential differences from the ones of $\mathbf J^1$ and $\mathbf J^2$. Thus we omit it.

As for \eqref{psia2}, as in \cite{seltes, huhoh}, we decompose $A_j, j = 1, 2$ to reveal the null structure of $A_{\mu}R^{\mu}_{\pm_{1}}\psi_{1,\pm_{1}}$.
\begin{equation*}
A_{j} = A_{j}^{df}+A_{j}^{cf},	
\end{equation*}
	where
	\begin{eqnarray*}
	A_{1}^{df}dx^{1}+A_{2}^{df}dx^{2} & = & (-\Delta)^{-1}(\partial_{2}({\rm curl}A)dx^{1}-\partial_{1}({\rm curl}A)dx^{2}),\\
	A_{1}^{cf}dx^{1}+A_{2}^{cf}dx^{2} & = & -(-\Delta)^{-1}(\partial_{1}({\rm div}A)dx^{1}+\partial_{2}({\rm div}A)dx^{2}).	
	\end{eqnarray*}
Then we obtain
$$
\sum_{\pm_1} A_{j}^{df}R_{\pm_{1}}^{j}\psi_{1, \pm_{1}} = \sum_{\pm_1, \pm_2}Q^{12}_{\pm_{2},\pm_{1}}(B_{\pm_{2}},\psi_{1, \pm_{1}}),
$$
where
$$
B_{\pm_2} = R_{\pm_2,1}A_{2,\pm_2}-R_{\pm_2,2}A_{1,\pm_2}.
$$
Also, we have
$$
\sum_{\pm_1}\left(A_{0}R^{0}_{\pm_{1}}\psi_{1, \pm_{1}} + A^{cf}_{j}R^{j}_{\pm_{1}}\psi_{1, \pm_{1}}\right) = \sum_{\pm_1, \pm_2} Q^{0}_{\pm_{1},\pm_{2}}(\psi_{1, \pm_{1}},A_{0,\pm_{2}}).
$$

Hence, to prove \eqref{psia2}, it suffices to show that
\begin{eqnarray*}
\|\Pi_{\pm_{0}}(Q^{12}_{\pm_{2},\pm_{1}}(B_{\pm_{2}},\psi_{1,\pm_{1}}))\|_{\mathcal{B}^{\frac{1}{4},-\frac{1}{2};1}_{\pm_{0}}} & \lesssim & \|\psi_{1,\pm_{1}}\|_{\mathcal{B}^{\frac{1}{4},\frac{1}{2};1}_{\pm_{1}}} \|B_{\pm_{2}}\|_{\mathcal{B}^{\frac{1}{4},\frac{1}{2};1}_{\pm_{2}}},	\\
\|\Pi_{\pm_{0}}(Q^{0}_{\pm_{1},\pm_{2}}(\psi_{1,\pm_{1}},A_{0,\pm_{2}}))\|_{\mathcal{B}^{\frac{1}{4},-\frac{1}{2};1}_{\pm_{0}}} & \lesssim & \|\psi_{1,\pm_{1}}\|_{\mathcal{B}^{\frac{1}{4},\frac{1}{2};1}_{\pm_{1}}} \|A_{0,\pm_{2}}\|_{\mathcal{B}^{\frac{1}{4},\frac{1}{2};1}_{\pm_{2}}}.
\end{eqnarray*}
The proofs are now straightforward from the previous argument. We just need to replace $\psi_{2,\pm_{2}}$ appearing in $\mathbf J^1$ and $\mathbf J^2$ by $B_{\pm_{2}}$. This completes the proof of Theorem \ref{lwp}.

\section{Failure of the smoothness}
\subsection{Proof of Theorem \ref{fail}} To show the failure of smoothness we adopt the argument of \cite{mst, hele}.
Let us consider the system of equations:
\begin{align}\label{csd-delta}
 \left\{
 \begin{array}{l}
 \Box A_{\nu}  =  \partial^{\mu}(-2\epsilon_{\mu\nu\lambda}\psi^{\dagger}\alpha^{\lambda}\psi),\\
 -i\alpha^{\mu}\partial_{\mu} \psi  =  A_{\mu}\alpha^{\mu}\psi, \\
 \psi(0) = \delta\psi_0,\quad A_\mu(0) = \delta a_\mu,\quad a_j = 0\\
 \partial_0 A_0(0) = 0, \quad \partial_0 A_j(0) = \delta\partial_j a_0 - 2\delta^2\epsilon_{0jk}\psi_0^\dagger\alpha^k\psi_0,
 \end{array}\right.\end{align}
where $0 < \delta \ll 1$. Let us denote the local solution of \eqref{csd-delta} by $(A_\nu(\delta, t), \psi(\delta, t))$. If the flow is $C^2$ at the origin in $H^s$, then it follows from \eqref{anu} and \eqref{psi} that
$$
\partial_\delta^2 A_\nu(0, t) = 2\sum_{\pm_1, \pm_2, \pm_3}\int_{0}^{t}U_{\pm_1}(t-t')\left[R_{\pm_1}^{\mu}(\epsilon_{\mu\nu\lambda}(U_{\pm_2}(t')\psi_0)^{\dagger}\alpha^{\lambda}U_{\pm_3}(t')\psi_0)\right]dt'
$$
and
$$
\partial_\delta^2 \psi(0, t) = 2i \sum_{\pm_1, \pm_2,\pm_3}\int_0^t U_{\pm_1}(t-t')\Pi_{\pm_1}(A_{\mu, \pm_2}^{\rm hom}(t') \alpha^\mu U_{\pm_3}(t')\psi_0)\,dt',
$$
where $U_\pm(t) := e^{-\pm it D}$, $A_{0, \pm}^{\rm hom} = \frac12 U_{\pm}(t) a_0$, and $A_{j, \pm}^{\rm hom} = -\pm\frac12U_{\pm}(t)\frac{\partial_j}{iD}a_0$.
From the $C^2$ smoothness we have that for a local existence time $T$
\begin{align}\begin{aligned}\label{smoothness}
&\sup_{0 \le t \le T}\left\|\sum_{\pm_1, \pm_2, \pm_3}\int_{0}^{t}U_{\pm_1}(t-t')\left[R_{\pm_1}^{\mu}(\epsilon_{\mu\nu\lambda}(U_{\pm_2}(t')\psi_0)^{\dagger}\alpha^{\lambda}U_{\pm_3}(t')\psi_0)\right]dt'\right\|_{H^s} \lesssim \|\psi_0\|_{H^s}^2,\\
&\sup_{0 \le t \le T}\left\|\sum_{\pm_1, \pm_2,\pm_3}\int_0^t U_{\pm_1}(t-t')\Pi_{\pm_1}(A_{\mu, \pm_2}^{\rm hom}(t') \alpha^\mu U_{\pm_3}(t')\psi_0)\,dt'\right\|_{H^s} \lesssim \|a_0\|_{H^s}\|\psi_0\|_{H^s}.
\end{aligned}\end{align}
However, we show that \eqref{smoothness} fails for $s < 0$.
\begin{proof}[Proof of failure of \eqref{smoothness}]
 Let $\psi_0 = \left(
                               \begin{array}{c}
                                 \frac1{\sqrt2}\phi \\
                                 -\frac1{\sqrt2}\phi \\
                               \end{array}
                             \right)$

and let
$$
F_\nu(t) := \sum_{\pm_1, \pm_2, \pm_3}\int_{0}^{t}U_{\pm_1}(t-t')\left[R_{\pm_1}^{\mu}(\epsilon_{\mu\nu\lambda}(U_{\pm_2}(t')\psi_0)^{\dagger}\alpha^{\lambda}U_{\pm_3}(t')\psi_0)\right]dt'.
$$
In particular,
\begin{align*}
\widehat{F_2}(t, \xi) &= \sum_{\pm_1, \pm_2, \pm_3}\int_{0}^{t}e^{-\pm_1i(t-t')|\xi|} \left(-1 + \pm_1\frac{\xi_1}{|\xi|}\right)\mathcal F_x(\overline{U_{\pm_2}(t')\phi}U_{\pm_3}(t')\phi)(\xi)\,dt'.
\end{align*}

Set $\phi = \mathcal F_x^{-1}(\chi_{W_\lambda})$ for $\lambda \gg 1$ and $W_\lambda = \{\xi = (\xi_1, \xi_2) : |\xi_1 - \lambda| \le \lambda^\frac12,\quad |\xi_2| \le \lambda^\frac12\}$.
Then
\begin{align*}
\widehat{F_2}(t, \xi) = -\sum_{\pm_1, \pm_2, \pm_3}e^{-\pm_1it|\xi|}\mathbf m_{123}(t, \xi, \eta) \left(1 - \pm_1\frac{\xi_1}{|\xi|}\right)\int_{\mathbb R^2} \chi_{W_\lambda}(\eta-\xi)\chi_{W_\lambda}(\eta)\,d\eta.
\end{align*}
where
\begin{align*}
\mathbf m_{123}(t, \xi, \eta) &:= \frac{e^{it\omega_{123}(\xi, \eta)}-1}{i\omega_{123}(\xi, \eta)},\quad \omega_{123}(\xi, \eta) = \pm_1|\xi| \pm_2|\xi-\eta| - \pm_3|\eta|.
\end{align*}
From the support condition we deduce that $|\xi| \lesssim \lambda^\frac12$. If $\pm_2 = \pm_3$, then $|\omega_{123}| \lesssim |\xi|$. Otherwise, $|\omega_{123}| \sim \lambda$.
Hence taking $t = \varepsilon \lambda^{-\frac12}$ for $0 < \varepsilon \ll 1$ and $\xi \in W_\lambda^* := \{3\xi_1^2 \le \xi_2^2\} \cap \{ |\xi| \lesssim \lambda^\frac12,\;\; |\xi_2| \sim \lambda^\frac12\} $, we get
$\mathbf m_{123}(t, \xi, \eta) = t(1 + O_{\pm}(\varepsilon))$
and $\sum_{\pm_1} e^{-\pm_1it|\xi|}(1 - \pm_1 \frac{\xi_1}{|\xi|}) = 2\cos(t|\xi|) + 2i\frac{\xi_1}{|\xi|}\sin(t|\xi|) = 2 + O(\varepsilon)$. Here $O_{\pm}(a)$ denotes that $O(a)$ depending on $\pm_j, j = 1, 2, 3$.
Thus we have
\begin{align*}
&\left|\sum_{\pm_1, \pm_2, \pm_3}e^{-\pm_1it|\xi|}\mathbf m_{123}(t, \xi, \eta) \left(1 - \pm_1\frac{\xi_1}{|\xi|}\right)\int_{\mathbb R^2} \chi_{W_\lambda}(\eta-\xi)\chi_{W_\lambda}(\eta)\,d\eta\,\right|\\
&\qquad \ge \left|\sum_{\pm_1, \pm_2 = \pm_3}e^{-\pm_1it|\xi|}\int_{\mathbb R^2}t(1 + O_{\pm}(\varepsilon))\left(1 - \pm_1\frac{\xi_1}{|\xi|}\right) \chi_{W_\lambda}(\eta-\xi)\chi_{W_\lambda}(\eta)\,d\eta\,\right|\\
&\qquad\qquad - C\sum_{\pm_1, \pm_2 \neq \pm_3}\lambda^{-1}\int_{\mathbb R^2} \chi_{W_\lambda}(\eta-\xi)\chi_{W_\lambda}(\eta)\,d\eta\\
&\qquad \ge 4t\int_{\mathbb R^2} \chi_{W_\lambda}(\eta-\xi)\chi_{W_\lambda}(\eta)\,d\eta - C(\varepsilon t+\lambda^{-1})\int_{\mathbb R^2} \chi_{W_\lambda}(\eta-\xi)\chi_{W_\lambda}(\eta)\,d\eta \\
&\qquad \gtrsim t\lambda \chi_{W_\lambda^*}(\xi) = \varepsilon \lambda^\frac12\chi_{W_\lambda^*}(\xi).
\end{align*}

Since $|(\widehat{F_0}(t, \xi), \widehat{F_1}(t, \xi), \widehat{F_2}(t, \xi))| \ge |\widehat F_2(t, \xi)|$, we finally get
$$
\varepsilon  \lambda^{\frac s2 + 1} \le \|F_2\|_{H^s} \le \|(F_0, F_1, F_2)\|_{H^s} \lesssim \|\psi_0\|_{H^s}^2 \lesssim \lambda^{2s + 1}.
$$
This means that the first part of \eqref{smoothness} fails for $s < 0$.

For the second part we take $a_0$
such that $\widehat{a_0} = \chi_{\widetilde W_\lambda}$ with $\widetilde W_\lambda = \{|\xi_1 - \lambda| \le \lambda^\frac14,\;\;|\xi_2|\le \lambda^\frac14\}$. Then by the same argument as above we have that for $\xi \in W_\lambda^{**} = \{|\xi_1 - 2\lambda| \le \frac1{10}\lambda^\frac12,\;\;|\xi_2| \le \frac1{10}\lambda^\frac12\}$, $t = \varepsilon \lambda^{-\frac12}$, and $0 < \varepsilon \ll 1$,
\begin{align*}
&\left|\sum_{\pm_1, \pm_2,\pm_3}\mathcal F_x\left(\int_0^t U_{\pm_1}(t-t')\Pi_{\pm_1}(A_{\mu, \pm_2}^{\rm hom}(t') \alpha^\mu U_{\pm_3}(t')\psi_0)\,dt'\right)(\xi)\right| \gtrsim \\
&\qquad\qquad \frac1{4\sqrt2}\left|\sum_{\pm_1, \pm_2 , \pm_3} e^{-\pm_1it|\xi|}\left(1 - \pm_1\frac{\xi_1}{|\xi|}\right)\int \widetilde{\mathbf m}_{123}(t, \xi, \eta)\left(1 \pm_2\frac{\xi_1- \eta_1}{|\xi-\eta|}\right) \chi_{\widetilde W_\lambda}(\xi-\eta)\chi_{W_\lambda}(\eta)\,d\eta\right|\\
&\qquad\qquad - \frac1{4\sqrt2}\left|\sum_{\pm_1, \pm_2 , \pm_3} e^{-\pm_1it|\xi|}\left( - \pm_1\frac{\xi_2}{|\xi|}\right)\int \widetilde{\mathbf m}_{123}(t, \xi, \eta)\left(\pm_2\frac{\xi_2- \eta_2}{|\xi-\eta|}\right) \chi_{\widetilde W_\lambda}(\xi-\eta)\chi_{W_\lambda}(\eta)\,d\eta\right|\\
&\qquad\qquad \gtrsim \varepsilon \lambda^\frac12\chi_{W_\lambda^{**}}(\xi),
\end{align*}
where $$\widetilde{\mathbf m}_{123} = \frac{e^{it(\pm_1|\xi| - \pm_2|\xi-\eta| - \pm_3|\eta|)}-1}{i(\pm_1|\xi| - \pm_2|\xi-\eta| - \pm_3|\eta|)}.$$
Taking $L_\xi^2((1+|\xi|)^{2s}d\xi)$ norm on both sides, we must have
$$
\varepsilon \lambda^{s+ 1} \lesssim \|a_0\|_{H^s}\|\psi_0\|_{H^s} \lesssim \lambda^{2s + \frac34},
$$
which implies the failure of the second part of \eqref{smoothness} for $s < \frac14$. This completes the proof of Theorem \ref{fail}.
\end{proof}

\subsection{ Proof of Theorem \ref{illp}}
In this subsection we show the smoothness failure of flow of \eqref{single} in $H^s$ for $s < 0$. The Cauchy problem \eqref{single} is equivalent to solving the integral equation
\begin{align*}
\psi(t) &= \sum_{\pm}\psi_\pm^{\rm hom} - \frac12\sum_{\pm_1, \pm_2}\int_0^t U_{\pm_1}(t-t')\Pi_{\pm_1} \left(\alpha^\nu  \left(U_{\pm_2}(t') (\pm_2 \frac1{D}\epsilon_{0\nu\lambda}\psi_0^\dagger \alpha^\lambda \psi_0)\right) \psi \right)\,dt'\\
&\qquad+ i\sum_{\pm_1, \pm_2}\int_0^t U_{\pm_1}(t-t')\Pi_{\pm_1}\left(\alpha^\nu \int_0^{t'} U_{\pm_2}(t'-t'')[R_{\pm_2}^\mu \epsilon_{\mu\nu\lambda}\psi^\dagger(t'') \alpha^\lambda \psi(t'')]\,dt'' \psi(t')   \right)\,dt',
\end{align*}
where $U_{\pm}(t) = e^{-\pm itD}$. In fact, it suffices to show the following.
\begin{prop}\label{ill-contraction}
Let $s < 0$. Then for a fixed $T > 0$ the inequality
\begin{align}\label{fail-cubic}
\sup_{t \in [0, T]}\left\| \mathcal N(\varphi)(t)\right\|_{H^s} \lesssim \|\varphi\|_{H^s}^3
\end{align}
fails to hold for any $\varphi \in H^s(\mathbb R^2)$, where $\mathcal N(\varphi)(t) = \sum_{\pm_j, j = 1,\cdots, 5} \mathcal N_{1\cdots5}$ and
$$
\mathcal N_{1\cdots5}(t) = \int_0^t U_{\pm_1}(t-t') \Pi_{\pm_1}\left(\alpha^\nu \int_0^{t'} U_{\pm_2}(t'-t'')\big[R_{\pm_2}^\mu \epsilon_{\mu\nu\lambda}(U_{\pm_3}(t'')\varphi)^\dagger\alpha^\lambda U_{\pm_4}(t'')\varphi\big]
\,dt'' U_{\pm_5}(t')\varphi \right )\,dt'.
$$
\end{prop}
\begin{rem}
By the same argument as described below one can also show that the estimate $\|\mathcal N_0(\varphi)\|_{H^s} \lesssim \|\varphi\|_{H^s}^3$ fails for $s < 0$. Here $\mathcal N_0(\varphi)$ is defined by
$$
\sum_{\pm_1, \pm_2, \pm_3}\int_0^t U_{\pm_1}(t-t')\Pi_{\pm_1} \left(\alpha^\nu  \left(U_{\pm_2}(t') (\pm_2 \frac1{D}\epsilon_{0\nu\lambda}\varphi^\dagger \alpha^\lambda \varphi)\right) U_{\pm_3}(t')\varphi \right)\,dt'.
$$
\end{rem}

\begin{proof} Let $\varphi = \left(
                               \begin{array}{c}
                                 \frac1{\sqrt2}\phi \\
                                 -\frac1{\sqrt2}\phi \\
                               \end{array}
                             \right)
$.
Then
\begin{align*}
\alpha^\nu  \big[R_{\pm_2}^\mu \epsilon_{\mu\nu\lambda}&(U_{\pm_3}(t'')\varphi)^\dagger\alpha^\lambda U_{\pm_4}(t'')\varphi\big]
\; U_{\pm_5}(t')\varphi \\
& = -\frac i{\sqrt 2} \overline{U_{\pm_3}(t'')\phi}\,U_{\pm_4}(t'')\phi U_{\pm_5}(t')\phi\left(
                               \begin{array}{c}
                                 1 \\
                                 1 \\
                               \end{array}
                             \right)\\
                              &\quad + \frac{1}{\sqrt2} R_{\pm_2}^1(\overline{U_{\pm_3}(t'')\phi}\,U_{\pm_4}(t'')\phi)U_{\pm_5}(t')\phi \left(
                               \begin{array}{c}
                                 1 \\
                                 -1 \\
                               \end{array}
                             \right).
\end{align*}
Since $\Pi_{\pm_1} = \frac12(I_{2\times 2} + R_{\pm_1, j} \alpha^j)$, we get
\begin{align*}
&\Pi_{\pm_1}\left(\alpha^\nu  U_{\pm_2}(t'-t'')\big[R_{\pm_2}^\mu \epsilon_{\mu\nu\lambda}(U_{\pm_3}(t'')\varphi)^\dagger\alpha^\lambda U_{\pm_4}(t'')\varphi\big]
\; U_{\pm_5}(t')\varphi \right)\\
&= -\frac i{2\sqrt 2} U_{\pm_2}(t'-t'')(\overline{U_{\pm_3}(t'')\phi}\,U_{\pm_4}(t''))\phi \;U_{\pm_5}(t')\phi\left(
                               \begin{array}{c}
                                 1 \\
                                 1 \\
                               \end{array}
                             \right)\\
&\qquad +  \frac{1}{2\sqrt2} U_{\pm_2}(t'-t'')[R_{\pm_2}^1(\overline{U_{\pm_3}(t'')\phi}\,U_{\pm_4}(t'')]\phi)\;U_{\pm_5}(t')\phi \left(
                               \begin{array}{c}
                                 1 \\
                                 -1 \\
                               \end{array}
                             \right)\\
&\qquad -\frac{i}{2\sqrt2} R_{\pm_1,1 }\left(U_{\pm_2}(t'-t'')(\overline{U_{\pm_3}(t'')\phi}\,U_{\pm_4}(t''))\phi\; U_{\pm_5}(t')\phi\right) \left(
                               \begin{array}{c}
                                1 \\
                                1 \\
                               \end{array}
                             \right)\\
&\qquad + \frac{1}{2\sqrt2} R_{\pm_1,1 }\left(U_{\pm_2}(t'-t'')[R_{\pm_2}^1(\overline{U_{\pm_3}(t'')\phi}\,U_{\pm_4}(t'')]\phi)\;U_{\pm_5}(t')\phi\right) \left(
                               \begin{array}{c}
                                -1 \\
                                1 \\
                               \end{array}
                             \right)\\
&\qquad - \frac{i}{2\sqrt2} R_{\pm_1,2 }\left(U_{\pm_2}(t'-t'')[\overline{U_{\pm_3}(t'')\phi}\,U_{\pm_4}(t'')\phi)]\;U_{\pm_5}(t')\phi\right) \left(
                               \begin{array}{c}
                                -1 \\
                                1 \\
                               \end{array}
                             \right)\\
&\qquad +  \frac{i}{2\sqrt2} R_{\pm_1,2 }\left(U_{\pm_2}(t'-t'')(R_{\pm_2}^1[\overline{U_{\pm_3}(t'')\phi}\,U_{\pm_4}(t'')\phi)])\;U_{\pm_5}(t')\phi\right) \left(
                               \begin{array}{c}
                                1 \\
                                1 \\
                               \end{array}
                             \right)\\
                             &=: -\frac i{2\sqrt2} \mathcal B^1 + \frac1{2\sqrt2}\mathcal B^2 -\frac i{2\sqrt2}\mathcal B^3 + \frac1{2\sqrt2}\mathcal B^4 - \frac i{2\sqrt2} \mathcal B^5 + \frac i{2\sqrt2} \mathcal B^6.
\end{align*}
Therefore
\begin{align*}
\mathcal N_{1\cdots5}(t) &= -\frac{i}{2\sqrt2}\int_0^t \!\!\!\int_0^{t'}  U_{\pm_1}(t-t')(\mathcal B^1 + \mathcal B^3 + \mathcal B^5 - \mathcal B^6)\,dt''dt'\\
&\qquad + \frac1{2\sqrt2}\int_0^t \!\!\!\int_0^{t'}  U_{\pm_1}(t-t')(\mathcal B^2 + \mathcal B^4)\,dt''dt'.
\end{align*}

Now let ${}^{\mathbf 1}\mathcal N_{1\cdots5}(t)$ be the first component of $\mathcal N_{1\cdots5}(t)$. Then the failure of \eqref{fail-cubic} is reduced to the one of the following:
\begin{align}\label{fail-n12345}
\left\|\sum_{\pm_j, j= 1,\cdots, 5} {}^{\mathbf 1} N_{1\cdots5}(t)\right\|_{H^s} \lesssim \|\phi\|_{H^s}^3.
\end{align}

Let $\lambda \gg 1$. Define $W_\lambda = \{\xi = (\xi_1, \xi_2) \in  \mathbb R^2 : | \xi_1 - \lambda| \le \lambda^\frac12, |\xi_2| \le \lambda^\frac12\}$.
Let $\phi = \mathcal F_\xi^{-1}(\chi_{W_\lambda})$. Then by taking space Fourier transform $\mathcal F_x$ we get
\begin{align}\begin{aligned}\label{m123}
&\mathcal F_x ({}^{\mathbf 1}\mathcal N_{1\cdots5}) (t, \xi)\\
& =  -\frac{i}{2\sqrt2}\left(1 \pm_1 \frac{\xi_1}{|\xi|}\right)\int\!\!\!\int_{\mathbb R^2 \times \mathbb R^2} \mathbf m_{1\cdots5}(t, \xi, \eta, \zeta)\chi_{W_\lambda}(\zeta-\eta)\chi_{W_\lambda}(\zeta)\chi_{W_\lambda}(\xi-\eta)\,d\eta d\zeta\\
&\qquad + \frac{i}{2\sqrt2}\left(\pm_1 \frac{\xi_2}{|\xi|}\right)\int\!\!\!\int_{\mathbb R^2 \times \mathbb R^2} \mathbf m_{1\cdots5}(t, \xi, \eta, \zeta)\left(1-\pm_2\frac{\eta_1}{|\eta|}\right) \chi_{W_\lambda}(\zeta-\eta) \chi_{W_\lambda}(\zeta) \chi_{W_\lambda}(\xi-\eta)\,d\eta d\zeta\\
&\qquad + \frac1{2\sqrt2}\left(1 - \pm_1 \frac{\xi_1}{|\xi|}\right) \int\!\!\!\int_{\mathbb R^2 \times \mathbb R^2}\mathbf m_{1\cdots5}(t, \xi, \eta, \zeta)\left(-\pm_2\frac{\eta_1}{|\eta|}\right)\chi_{W_\lambda}(\zeta-\eta)\chi_{W_\lambda}(\zeta)\chi_{W_\lambda}(\xi-\eta)\,d\eta d\zeta\\
&=: \mathcal M_{1\cdots5}^1 + \mathcal M_{1\cdots5}^2 + \mathcal M_{1\cdots5}^3,
\end{aligned}\end{align}
where
\begin{align*}
\mathbf m_{1\cdots5}(t, \xi, \eta, \zeta) &:= \int_0^t \!\!\!\int_0^{t'} e^{-\pm_1i(t-t')|\xi| -\pm_2i(t'-t'')|\eta|+ \pm_3it''|\zeta-\eta| - \pm_4it''|\zeta| - \pm_5it'|\xi-\eta|}\,dt''dt'\\
&= \frac{e^{-\pm_1it|\xi|}}{i\omega_0}\left(\frac{e^{it\omega_1}-1}{i\omega_1} - \frac{e^{it\omega_2}-1}{i\omega_2}\right)
\end{align*}
and
\begin{align*}
\omega_0 &= \pm_2|\eta| \pm_3 |\eta - \zeta| - \pm_4|\zeta|,\\
\omega_1 &= \pm_1|\xi| - \pm_5|\xi-\eta| \pm_3 |\eta - \zeta| - \pm_4|\zeta|,\\
\omega_2 &= \pm_1|\xi| - \pm_2 |\eta| - \pm_5|\xi-\eta|.
\end{align*}

From the support condition it follows that $\eta \in \{\eta : |\eta| \lesssim \lambda^\frac12\}$, provided $\lambda \gg 1$ and $\xi \in W_\lambda$.
If $\pm_1 = \pm_5 = \pm_3 = \pm_4$, then $|\omega_1| \lesssim 1$. If $\pm_1 =\pm_5$ and $\pm_3=\pm_4$, $\pm_1 \neq \pm_3$, then $|\omega_0|\lesssim |\eta|, |\omega_1|\lesssim |\eta|,|\omega_2|\lesssim |\eta|$. If $\pm_1=\pm_5=\pm_3$ and $\pm_3\neq\pm_4$, or $\pm_1=\pm_5=\pm_4$ and $\pm_4\neq\pm_3$, then $|\omega_0|\sim\lambda, |\omega_1|\sim\lambda$ and $|\omega_2|\lesssim |\eta|$.\\

Let us set $t = \varepsilon \lambda^{-\frac12}$ and $\lambda = 4k^2\pi^2\varepsilon^{-2}$ for $1 \ll k \in \mathbb N$. Then we first have
\begin{align*}
\sum_{\pm_1 = \pm_5, \pm_3 = \pm_4}&\left(1 \pm_1 \frac{\xi_1}{|\xi|}\right)\mathbf m_{1\cdots5}(t, \xi, \eta, \zeta)\\
 &=  \sum_{\pm_1 = \pm_5, \pm_3 = \pm_4}\left(1 \pm_1 \frac{\xi_1}{|\xi|}\right)\frac{t^2e^{-\pm_1it|\xi|}}{2}(1 + O_\pm(\varepsilon)) \\
 &= \sum_{\pm_1 = \pm_5, \pm_3 = \pm_4}\left(1 \pm_1 \frac{\xi_1}{|\xi|}\right)\frac{t^2e^{-\pm_1it|\xi|}}{2} + \sum_{\pm_1 = \pm_5, \pm_3 = \pm_4}\left(1 \pm_1 \frac{\xi_1}{|\xi|}\right)\frac{t^2e^{-\pm_1it|\xi|}}{2}\;O_\pm(\varepsilon)\\
&= 4t^2\left(\cos(t|\xi|) - i\frac{\xi_1}{|\xi|}\sin(t|\xi|)\right) + \sum_{\pm_1 = \pm_5, \pm_3 = \pm_4}\left(1 \pm_1 \frac{\xi_1}{|\xi|}\right)\frac{t^2e^{-\pm_1it|\xi|}}{2}\;O_\pm(\varepsilon)\\
&= 4t^2(1 + O(\varepsilon)) + \sum_{\pm_1 = \pm_5, \pm_3 = \pm_4}\left(1 \pm_1 \frac{\xi_1}{|\xi|}\right)\frac{t^2e^{-\pm_1it|\xi|}}{2}\;O_\pm(\varepsilon).
\end{align*}
Here $O_\pm(\varepsilon)$ denotes $O(\varepsilon)$ depending on $\pm_j, j= 1, \cdots, 5$.
Hence we get
\begin{align*}
\sum_{\pm_1 = \pm_5, \pm_3 = \pm_4}\mathcal M_{1\cdots5}^1 =  -\frac{2it^2}{\sqrt2} \int\!\!\!\int_{\mathbb R^2 \times \mathbb R^2}\chi_{W_\lambda}(\zeta-\eta)\chi_{W_\lambda}(\zeta)\chi_{W_\lambda}(\xi-\eta)\,d\eta d\zeta + O(\varepsilon^3\lambda).
\end{align*}
This gives us
\begin{align}\label{m1}
\left|\sum_{\pm_1 = \pm_5, \pm_3 = \pm_4}\mathcal M_{1\cdots5}^1\right| \gtrsim \varepsilon^2\lambda.
\end{align}

Similarly, we have
\begin{align}\begin{aligned}\label{m2}
&\left|\sum_{\pm_1 = \pm_5, \pm_3 = \pm_4}\mathcal M_{1\cdots5}^2\right|\\
&\qquad =  \left|-\frac{2t^2}{\sqrt2}\frac{\xi_2}{|\xi|}\sin(t|\xi|) \int\!\!\!\int_{\mathbb R^2 \times \mathbb R^2}\chi_{W_\lambda}(\zeta-\eta)\chi_{W_\lambda}(\zeta)\chi_{W_\lambda}(\xi-\eta)\,d\eta d\zeta\right| + O(\varepsilon^3\lambda) \lesssim \varepsilon^3\lambda.
\end{aligned}\end{align}

For $\mathcal{M}^{3}_{1\cdots5}$,
\begin{align*}
\sum_{\pm_1 = \pm_5, \pm_3 = \pm_4}&\left(1 -\pm_1 \frac{\xi_1}{|\xi|}\right)\mathbf m_{1\cdots5}(t, \xi, \eta, \zeta)\left(-\pm_{2}\frac{\eta_{1}}{|\eta|}\right)\\
 &=  \sum_{\pm_1 = \pm_5, \pm_3 = \pm_4}\left(1 \pm_1 \frac{\xi_1}{|\xi|}\right)\frac{t^2e^{-\pm_1it|\xi|}}{2}(1 + O_\pm(\varepsilon))\left(-\pm_{2}\frac{\eta_{1}}{|\eta|}\right) \\
 &= \sum_{\pm_1 = \pm_5, \pm_3 = \pm_4}\left(1 \pm_1 \frac{\xi_1}{|\xi|}\right)\frac{t^2e^{-\pm_1it|\xi|}}{2}\;O_\pm(\varepsilon)\left(-\pm_{2}\frac{\eta_{1}}{|\eta|}\right), \quad\bigg(\because \sum_{\pm_{2}}\left(-\pm_{2}\frac{\eta_{1}}{|\eta|}\right)=0\bigg)
\end{align*}
and
$$
\left|\sum_{\pm_1 = \pm_5, \pm_3 = \pm_4}\left(1 \pm_1 \frac{\xi_1}{|\xi|}\right)\frac{t^2e^{-\pm_1it|\xi|}}{2}\;O_\pm(\varepsilon)\left(-\pm_{2}\frac{\eta_{1}}{|\eta|}\right) \right|\lesssim \varepsilon^{3}\lambda^{-1}.
$$
Thus
\begin{align}\label{m3}
\left|\sum_{\pm_1 = \pm_5, \pm_3 = \pm_4}\mathcal M_{1\cdots5}^2\right|  \lesssim \varepsilon^3\lambda.
\end{align}

In the case $\pm_1 = \pm_5$ and $\pm_3 \neq \pm_4$ $|\mathbf m_{1\cdots5}| \lesssim \lambda^{-1}(\lambda^{-1} + t)$ and hence
\begin{align}\label{m4}
\left|\sum_{\pm_1 = \pm_5, \pm_3 \neq \pm_4}\sum_{j = 1, 2, 3}\mathcal M_{1\cdots5}^j\right| \lesssim \varepsilon \lambda^\frac12.
\end{align}

On the other hand, we write $$\mathbf m_{1 \cdots 5} = -\frac{e^{-\pm_1it|\xi|}}{\omega_1\omega_2}\left(1 - e^{it\omega_1} + \omega_1e^{it\omega_2}\frac{e^{it\omega_0}-1}{\omega_0}\right).$$
We consider the cases $\pm_{1}\neq\pm_{5}$ with $\pm_{1}=\pm_{4},\pm_{5}=\pm_{3}$ and $\pm_{1}\neq\pm_{5}$ with $\pm_{5}=\pm_3=\pm_4$. Since $|1-e^{it\omega_1}|\le t\omega_1$, in the former case ($|\omega_0|\sim \lambda,|\omega_1|\lesssim\lambda,|\omega_2|\sim\lambda$), we get
$$\left|\mathbf{m}_{1\cdots5}\right| \le \left|\frac{e^{-\pm_1it|\xi|}}{\omega_1\omega_2} \right| \left(\left|t\omega_1\right|+\left|\omega_1 \frac{e^{it\omega_0}-1}{\omega_0}\right|\right)\lesssim \lambda^{-1}\left(t+\lambda^{-1}\right) \lesssim t\lambda^{-1}.$$
In the latter case ($|\omega_0|\lesssim|\eta|,|\omega_1|\sim\lambda,|\omega_2|\sim\lambda$), we get
$$
\left|\mathbf{m}_{1\cdots5}\right| \lesssim \frac{t}{|\omega_2|}\lesssim t\lambda^{-1}.
$$
Therefore, we estimate
\begin{align}\label{m5}
\left|\sum_{\pm_1 \neq \pm_5}\sum_{j = 1, 2, 3}\mathcal M_{1\cdots5}^j\right| \lesssim \varepsilon \lambda^\frac12.
\end{align}

Putting all the estimates \eqref{m1}--\eqref{m5} into \eqref{m123}, we finally get
$$
\bigg\|\sum_{\pm_j, j= 1,\cdots, 5} {}^{\mathbf 1} N_{1\cdots5}(t)\bigg\|_{H^s} \gtrsim \varepsilon^2\lambda^{\frac32 + s},
$$
provided $\varepsilon\lambda^\frac12 = 2\pi k \gg 1$ and $\varepsilon \ll 1$.
This gives us the failure of \eqref{fail-n12345} for $s < 0$ and thus \eqref{fail-cubic}. This completes the proof of Proposition \ref{ill-contraction}.

\end{proof}

\section{Appendix}

\subsection{Proof of Lemma \ref{energy-est}}

Let $G\in\mathcal{B}^{s,-\frac{1}{2};1}_{\pm}$ be an arbitrarily chosen representative of $F\in\mathcal{B}^{s,-\frac{1}{2};1}_{\pm}(S_{T})$. Then we need to show that $\|u\|_{\mathcal{B}^{s,\frac{1}{2};1}_{\pm}(S_{T})} \lesssim \|f\|_{B^{s}_{2,1}} + \|G\|_{\mathcal{B}^{s,-\frac{1}{2};1}_{\pm}(S_T)}$.
By density, we may assume $G\in\mathcal{S}(\mathbb{R}^{1+2})$.
Since $\mathcal{F}\lbrack e^{\mp itD}f\rbrack=\delta(\tau\pm|\xi|)\widehat{f}(\xi)$, by taking a $C_0^\infty(\mathbb R)$ function $\rho$ with value 1 on $(-T,T)$ we have
\begin{eqnarray*}
\|e^{\mp itD}f\|_{\mathcal{B}^{s,\frac{1}{2};1}_{\pm}(S_{T})}	 & \le & \|\rho e^{\mp itD}f\|_{\mathcal{B}^{s,\frac{1}{2};1}_{\pm}} = \sum_{N,L\ge1}N^{s}L^{1/2}\|\chi_{|\xi|\sim N}\chi_{|\tau\pm|\xi||\sim L}\widehat{\rho}(\tau\pm|\xi|)\widehat{f}(\xi)\|_{L^{2}_{\tau,\xi}}\\
& \le & \sum_{N,L\ge1}N^{s}L^{1/2}\|\chi_{|\cdot|\sim L}\widehat{\rho}\|_{L^{2}_{\tau}}\|\chi_{|\cdot|\sim N}\widehat{f}\|_{L^{2}_{\xi}} \le \|\rho\|_{B^{1/2}_{2,1}}\|f\|_{B^{s}_{2,1}} \lesssim \|f\|_{B^{s}_{2,1}}.
\end{eqnarray*}
Next, we write $w(t)=\int_{0}^{t}e^{\mp i(t-t')|\xi|}G(t')dt'$ and by taking Fourier transform in space we get
\begin{eqnarray*}
\widehat{w}(t,\xi) = \int_{0}^{t}e^{\mp i(t-t')|\xi|}\widehat{G}(t',\xi)dt' = \frac1{2\pi}\int_{0}^{t}e^{\mp i(t-t')|\xi|}\int e^{it'\lambda}\widetilde{G}(\lambda,\xi)d\lambda dt' = \frac1{2\pi}\int \frac{e^{it\lambda}-e^{\mp it|\xi|}}{i(\lambda\pm|\xi|)}\widetilde{G}(\lambda,\xi)d\lambda.
\end{eqnarray*}
Thus we have
$$
\widetilde{w}(\tau,\xi)=\frac{\widetilde{G}(\tau,\xi)}{i(\tau\pm|\xi|)}-\delta(\tau\pm|\xi|)\int\frac{\widetilde{G}(\lambda,\xi)}{i(\lambda\pm|\xi|)}d\lambda.
$$
Now we split $G=G_{1}+G_{2}$ corresponding to the Fourier domains $|\tau\pm|\xi||\lesssim 1$ and $|\tau\pm|\xi||\gg1$, respectively. We write $w=w_{1}+w_{2}$ accordingly. By Taylor's expansion, we get
$$\widehat{w_{1}}(t,\xi)=e^{\mp it|\xi|}\sum_{n\ge1}\int\frac{(it(\lambda\pm|\xi|))^{n}}{n!i(\lambda\pm|\xi|)}\chi_{|\tau\pm|\xi||\lesssim1}\widetilde{G}(\lambda,\xi)d\lambda.$$
Hence we write $w_{1}(t)=\sum_{n\ge1}\frac{t^{n}}{n!}e^{\mp itD}f_{n}$, where $\widehat{f_{n}}(\xi)=\int(i(\lambda\pm|\xi|))^{n-1}\chi_{|\lambda\pm|\xi||\lesssim1}\widetilde{G}(\lambda,\xi)d\lambda$.
Now from the support condition $\{|\lambda\pm|\xi||\lesssim1\}$ it is easy matter to see that $\|f_{n}\|_{B^{s}_{2,1}}\lesssim \|G\|_{\mathcal{B}^{s,-\frac{1}{2};1}_{\pm}}$. %Indeed,
%\begin{eqnarray*}
%||f_{n}||_{B^{s}_{2,1}} & = & \sum_{N\ge1}N^{s}||\chi_{|\xi|\sim N}\widehat{f_{n}}||_{L^{2}_{\xi}}\\
%& \lesssim & \sum_{N,L\ge1}N^{s}\bigg\lvert\bigg\lvert
%\chi_{|\xi|\sim N}\chi_{|\lambda\pm|\xi||\sim L}\int(i(\lambda\pm|\xi|))^{n-1}\chi_{|\lambda\pm|\xi||\lesssim1}\widetilde{G}(\lambda,\xi)d\lambda\bigg\rvert\bigg\rvert_{L^{2}_{\xi}}\\
%& \lesssim & \sum_{N\ge1}\sum_{L\lesssim1}N^{s}L^{-1}\bigg|\bigg|\chi_{|\xi|\sim N}\int\chi_{|\lambda\pm|\xi||\sim L}\widetilde{G}(\lambda,\xi)d\lambda\bigg|\bigg|_{L^{2}_{\xi}}\\
%& \lesssim & \sum_{N,L\ge1}N^{s}L^{-1}L^{1/2}||\chi_{|\xi|\sim N}\chi_{|\tau\pm|\xi||\sim L}\widetilde{G}||_{L^{2}_{\tau,\xi}} \\
%& = & ||G||_{\mathcal{B}^{s,-\frac{1}{2};1}_{\pm}}
%\end{eqnarray*}
Thus we get
\begin{eqnarray*}
\|w_{1}\|_{\mathcal{B}^{s,\frac{1}{2};1}_{\pm}(S_{T})} & \le &  \|\rho w_{1}\|_{\mathcal{B}^{s,\frac{1}{2};1}_{\pm}} \le \sum_{n\ge1}\frac{1}{n!}\|t^{n}\rho(t)e^{\mp itD}f_{n}\|_{\mathcal{B}^{s,\frac{1}{2};1}_{\pm}} \\
& \le & \sum_{n\ge1}\frac{1}{n!}\|t^{n}\rho(t)\|_{B^{1/2}_{2,1}}\|f_{n}\|_{B^{s}_{2,1}} \lesssim \bigg(\sum_{n\ge1}\frac{n2^{n-1}}{n!}\bigg)\|G\|_{\mathcal{B}^{s,-\frac{1}{2};1}_{\pm}} \lesssim \|G\|_{\mathcal{B}^{s,-\frac{1}{2};1}_{\pm}}.
\end{eqnarray*}
where we used $T \le 1$ and $\|t^{n}\rho(t)\|_{B^{1/2}_{2,1}}\lesssim \|t^{n}\rho(t)\|_{H^{1}}\lesssim n2^{n-1}$.

To estimate $w_{2}$, we split $w_{2}=a-b$ where
\begin{eqnarray*}
\widetilde{a}(\tau,\xi) & = & \frac{\chi_{|\tau\pm|\xi||\gg1}\widetilde{G}(\tau,\xi)}{i(\tau\pm|\xi|)}, \\
\widetilde{b}(\tau,\xi) & = & \delta(\tau\pm|\xi|)\int\frac{\chi_{|\lambda\pm|\xi||\gg1}\widetilde{G}(\lambda,\xi)}{i(\lambda\pm|\xi|)}d\lambda.
\end{eqnarray*}
Thus
\begin{eqnarray*}
\|a\|_{\mathcal{B}^{s,\frac{1}{2};1}_{\pm}(S_T)}  &\le& \|\rho a\|_{\mathcal{B}^{s,\frac{1}{2};1}_{\pm}} \lesssim  \sum_{N\ge1}N^{s}\sum_{L\gg1}L^{1/2}L^{-1}\|\chi_{|\xi|\sim N}\chi_{|\tau\pm|\xi||\sim L}\widetilde{G}\|_{L^{2}_{\tau,\xi}} \le \|G\|_{\mathcal{B}^{s,-\frac{1}{2};1}_{\pm}},\\
\|b\|_{\mathcal{B}^{s,\frac{1}{2};1}_{\pm}(S_{T})} &\le& \|\rho b\|_{\mathcal{B}^{s,\frac{1}{2};1}_{\pm}} \lesssim  \sum_{N\ge1}N^{s}\sum_{L\gg1}L^{-1}L^{1/2}\|\chi_{|\xi|\sim N}\chi_{|\tau\pm|\xi||\sim L}\widetilde{G}\|_{L^{2}_{\tau,\xi}} \le \|G\|_{\mathcal{B}^{s,-\frac{1}{2};1}_{\pm}}.	
\end{eqnarray*}
This completes the proof of Lemma \ref{energy-est}.

%%%%%%%%%%%%%%%%%%%%%%%%%%%%%%%%%%%%%%%%%%%%%%%%%%%%%%%%%%%%%%%%%%%%%%%%%%%%%%%%%%%%%%%%%%%%%%%%%%%%%%%%%%%%%%%%%%%%%%%%%%%%%%%%%%%%%%%%%%%%%
%%%%%%%%%%%%%%%%%%%%%%%%%%%%%%%%%%%%%%%%%%%%%%%%%%%%%%%%%%%%%%%%%%%%%%%%%%%%%%%%%%%%%%%%%%%%%%%%%%%%%%%%%%%%%%%%%%%%%%%%%%%%%%%%%%%%%%%%%%%%%
\section*{Acknowledgements}
This work was supported by NRF-2018R1D1A3B07047782(Republic of Korea).
%%%%%%%%%%%%%%%%%%%%%%%%%%%%%%%%%%%%%%%%%%%%%%%%%%%%%%%%%%%%%%%%%%%%%%%%%%%%%%%%%%%%%%%%%%%%%%%%%%%%%%%%%%%%%%%%%%%%%%%%%%%%%%%%%%%%%%%%%%%%%
%%%%%%%%%%%%%%%%%%%%%%%%%%%%%%%%%%%%%%%%%%%%%%%%%%%%%%%%%%%%%%%%%%%%%%%%%%%%%%%%%%%%%%%%%%%%%%%%%%%%%%%%%%%%%%%%%%%%%%%%%%%%%%%%%%%%%%%%%%%%%

%%%%%%%%%%%%%%%%%%%%%%%%%%%%%%%%%%%%%%%%%%%%%%%%%%%%%%%%%%%%%%%%%%%%%%%%%%%%%%%%%%%%%%%%%%%%%%%%%%%%%%%%%%%%%%%%%%%%%%%%%%%%%%%%%%%%%%%%%%%%%%%%%%%%%%%%%%%%%%%%%%%%%%%%%%%%%%%%%%%%%%%%%%%%%%%%%%%%%%%%%%%%%%%%%%%%%%%%%%%%%%%%%%%%%%%%%%%%%%%%%%%%%%%%%%%%%%%%%%%%%%%%%%%%%%%%%%%%%%%%%%%%%%%%%%%%%%%%%%%%%%%%%%%%

\end{document}